\documentclass[11pt]{article}
\usepackage{dsfont} 
\usepackage{amssymb} 
\usepackage{graphicx}
\usepackage{extarrows}

\usepackage{mathrsfs} 

\begin{document}
\title{{\bf\Large Effective dynamics of stochastic wave equation
with a random dynamical boundary condition
\thanks{This work was supported by the National Science
Foundation of China (Grants No. 10901115, 10971225 and 11071177) and the
Scientific Research Found of Science and Technology Bureau of
Sichuan Province (Grant No. 2012JQ0041 and 2010JY0057).  } } }
\author{Guanggan Chen\\
College of Mathematics and Software Science,\\
Sichuan Normal University,
Chengdu, 610068, China\\
\emph{E-mail: chenguanggan@hotmail.com}\\
\\
Jinqiao Duan\\
Institute for Pure and Applied Mathematics,\\
University of California,
Los Angeles, CA 90095, USA\\
\emph{E-mail: jduan@ipam.ucla.edu }\\ \& \\
Department of Applied Mathematics, \\
Illinois Institute of Technology,
Chicago, IL 60616, USA \\
\emph{E-mail: duan@iit.edu }\\
\\
Jian Zhang\\
College of Mathematics and Software Science, \\
Sichuan Normal University,
Chengdu, 610068, China\\
\emph{E-mail: zhangjiancdv@sina.com}
 }

\date{\today}
\maketitle
\baselineskip 6.8mm
\begin{center}
\begin{minipage}{120mm}
{\bf \large Abstract:} {\small This work is devoted to the effective
macroscopic dynamics of a   weakly damped stochastic
nonlinear wave equation with a random dynamical boundary condition.
The white noises are taken into account not only in the model
equation defined on a domain perforated with small holes, but also
in the dynamical boundary condition on the boundaries of the small
holes. An effective homogenized, macroscopic model is
derived in the sense of probability distribution, which is a new
stochastic wave equation on a unified domain, without small holes,
with a usual static boundary condition. }
\end{minipage}
\end{center}

\renewcommand{\theequation}{\thesection.\arabic{equation}}
\setcounter{equation}{0}

\vspace{0.5cm}

\section{ Introduction }

\quad\quad Nonlinear wave equations, as a class of important
mathematical models, describe the propagation of waves in certain
systems or media, such as sonic booms,   traffic flows,
  optic devices and quantum fields   (\cite{RS, Whitham}). In
the deterministic case, they have been studied extensively   due to
their wide applications in engineering and science (e.g., \cite{GV, L,
PS, S}). On bounded domains or
media, the effect of the boundary   often needs to be considered.
 Dirichlet, Neumann and Robin boundary conditions  are
called {\it static boundary conditions}, as they are not
involved with time derivatives of the system state variables. On the
contrary,   {\it dynamic boundary conditions} contain time derivatives of the system
state variables and arise  in many
physical problems   (see \cite{FG, GS, PR}).
\par
In some physical problems, such as wave propagation through the
atmosphere or the ocean, due to   stochastic force, uncertain parameters, random
sources and random boundary conditions,  the realistic models
take the random fluctuation into account
\cite{Chen1, Chow1, DZ1}. This leads
  to  stochastic nonlinear wave equations, which have drawn
quite attentions recently   \cite{Chen2, Chow1, Chow2, FW, LS07, LW, M,
ZYO}.
\par
In this paper, we are concerned with the effective, macroscopic
dynamics of the following ``microscopic" weakly damped stochastic
nonlinear wave equation with a random dynamical boundary condition
on a domain $D$ perforated with small holes
\begin{equation}\label{Eq1}
\left\{
\begin{array}{ll}
u^\varepsilon_{tt}+u^\varepsilon_t-\bigtriangleup
u^\varepsilon+u^\varepsilon-f(u^\varepsilon)=\dot{W_1} &\quad
\hbox{in}\; D^\varepsilon\times
[0, \tau^*),\\
\varepsilon^2\delta^\varepsilon_{tt}+\delta^\varepsilon_t+\varepsilon^2\delta^\varepsilon
=-\varepsilon^2u^\varepsilon_t+\varepsilon^2\dot{W_2} &\quad
\hbox{on}\;
\partial S^\varepsilon\times [0, \tau^*),\\
u^\varepsilon=0 &\quad \hbox{on}\; \partial D\times[0,
\tau^*),\\
\delta^\varepsilon_t=\frac{\partial u^\varepsilon}{\partial {\bf n}}
&\quad \hbox{on}\; \partial S^\varepsilon\times [0, \tau^*).
\end{array}
\right.
\end{equation}
Here $\varepsilon $ is a small positive parameter, and the domain $D^\varepsilon$ is a subset of an open bounded
domain $D$ in $\mathbb{R}^3$, obtained by removing $S^\varepsilon$,
the collection of small holes of size $\varepsilon$, periodically
distributed in $D$. Also, $W_1$ and $W_2$ are two independent Wiener processes. This will be given in details in the next
section.  The symbol $\tau^*$ denotes a stopping time on
$(0,+\infty)$, and $\frac{\partial }{\partial {\bf n}}$ denotes the
unit outer normal derivative on the boundary $\partial
S^\varepsilon$. In particular, in this paper we will only concern
with the case of the nonlinear term $f(u^\varepsilon)=\sin
u^\varepsilon$ (the Sine-Gordon equation).
\par
The system (\ref{Eq1}), when the white noises, $\dot{W}_1,
\dot{W}_2$, and the parameter $\varepsilon$ are absent, arises in
the modeling of gas dynamics in an open bounded domain $D$, with
points on boundary acting like a spring reacting to the excess
pressure of the gas (see \cite{Frigeri, MI}). In this deterministic
case, Beale \cite{Beale1, Beale2} and Mugnolo \cite{Mugnolo}
established the well-posedness and analyzed some properties of the
spectrum in some special cases. Cousin, Frota and Larkin \cite{CFL}
studied the global solvability and asymptotic behavior. Frigeri
\cite{Frigeri} considered   large time dynamical behavior.
Furthermore, for the stochastic system (\ref{Eq1}), when the
parameter $\varepsilon$ are absent, Chen and Zhang \cite{Chen2}
investigated  the long time behavior of the solutions.
\par

Homogenization  plays an
important role in understanding multiscale   phenomena
in material science, climate dynamics, chemistry and biology
\cite{CD, Tartar}. For the deterministic system defined on
heterogeneous media, there have been some relevant works  for   heat conduction \cite{NR1, NR2, TTM} and for wave propagation
  \cite{CDMZ, Timofte}.  Several authors also considered
 homogenization problems for the random partial differential
equations (PDEs with random coefficients) \cite{KP, PP} and for
the partial differential equations on randomly heterogeneous
domains \cite{BM, Z1, Z2}. However, for the stochastic partial differential equations
(PDEs with white noises), especially for the stochastic partial
differential equations with random dynamical boundary conditions, due
to the effect of both nonlinear dynamical boundary condition and the
nonclassical fluctuation of driving white noises, the study of
stochastic homogenization problem is still in its infancy (see
\cite{WCD, WD}).
\par
Therefore, in this paper, we are especially interested in the
stochastic homogenization problem of Equation (\ref{Eq1}). Our aim
is to establish the effective macroscopic equation of Equation
(\ref{Eq1}). For this purpose,   the key step is to verify the
compactness of the solutions in some function space for the
deterministic systems. But it does not hold for stochastic Equation
(\ref{Eq1}). Therefore, we will instead consider the tightness of
the distributions of the solutions,  so that the effective
macroscopic equation is established in the sense of probability
distribution. More precisely, we first analyze the microscopic model
Equation (\ref{Eq1}) to establish the well-posedness.   Since the
energy relation of this stochastic system does not directly imply
the a priori estimate of the solutions, we then introduce a pseudo
energy argument to infer almost sure boundedness of the solutions.
Furthermore, we use the a priori estimate to establish the tightness
of distribution of the solutions. Finally, we derive the effective
homogenized
  equation in the sense of probability distribution, which
is a new stochastic wave equation on a unified domain without small
holes but with a static boundary condition. The solutions of the
original model Equation (\ref{Eq1}) converge to those of the
effective homogenized   equation in probability
distribution, as the size of small holes $\varepsilon$ diminishes to
zero.
\par
This paper is organized as follows. In the next section, we will
formulate the basic setup of the homogenization problem. In
section 3, we will prove the well-posedness, almost sure boundedness
and tightness of distribution of the solutions for the microscopic
model Equation (\ref{Eq1}). In     section 4, we will derive
the effective homogenized   equation in probability
distribution.
\par

\renewcommand{\theequation}{\thesection.\arabic{equation}}
\setcounter{equation}{0}

\section{Basic setup of the problem}

\begin{figure}[ht]
\begin{center}
\includegraphics[width=14cm,height=6cm]{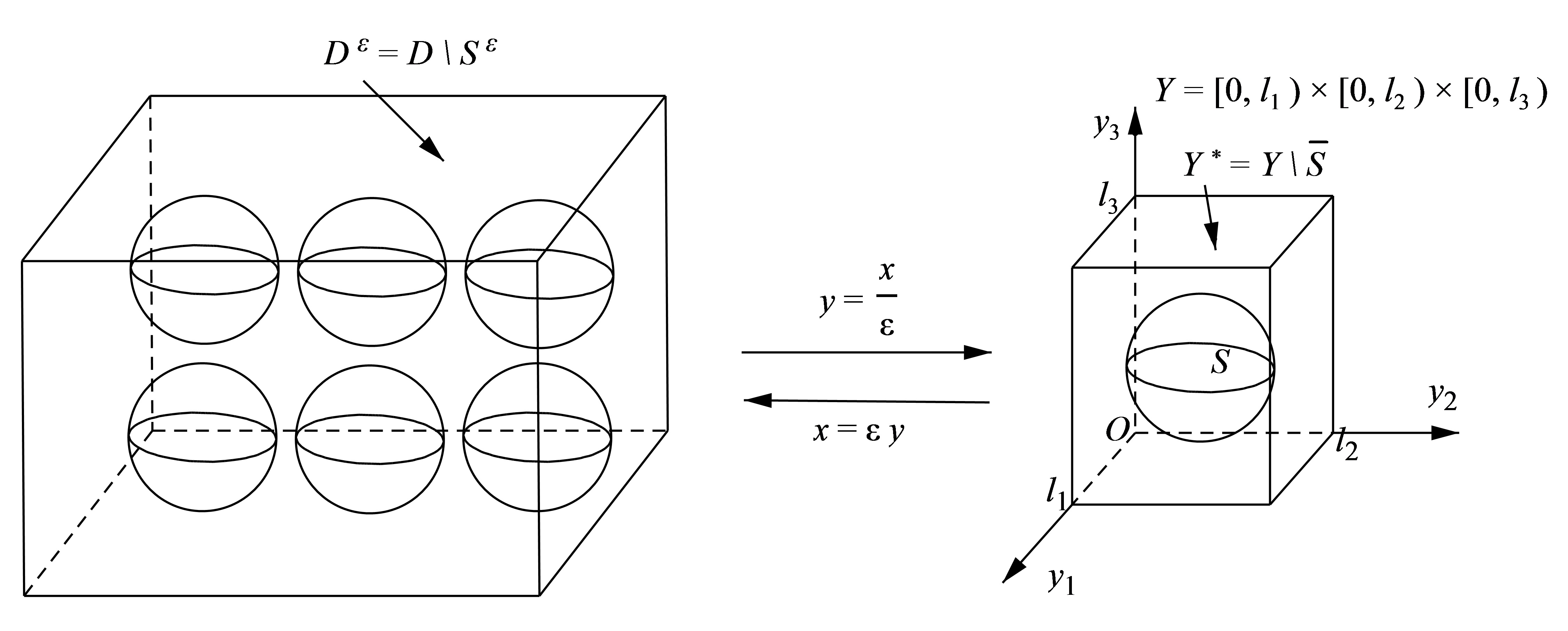}\vspace{-0.5cm}
\caption{ Geometric setup in $\mathbb{R}^3$}\label{f1}
\vspace{-0.5cm}
\end{center}
\end{figure}

\quad \quad Let the physical medium $D$ be an open bounded domain in
$\mathbb{R}^3$ with piece-wise smooth boundary $\partial D$, and let
$\varepsilon \in (0,1)$ be a small real parameter. Denote by $Y=[0,
l_1)\times [0, l_2)\times [0, l_3)$   a representative elementary
cell in $\mathbb{R}^3$ and let $S$ be an open subset of $Y$ with smooth
boundary $\partial S$ such that $S\subset Y$. The elementary cell
$Y$ and the small cavity of hole $S$ inside it are used to model
small scale obstacles or heterogeneities in a physical medium $D$.
Define $\varepsilon S=\{\varepsilon y: y\in S\}$ and
$S_{\varepsilon, {\bf k}}={{\bf k}l+\varepsilon S}$ with ${\bf
k}l=(k_1l_1, k_2l_2, k_3l_3)$ and ${\bf k}=(k_1, k_2, k_3)\in
\mathbb{Z}^3$. Let $S^\varepsilon$ be the set of all the holes
contained in $D$, i.e.,
$$
S^\varepsilon=\bigcup\{S_{\varepsilon, {\bf
k}}|\;\overline{S_{\varepsilon, {\bf k}}} \subset D, \hbox{and}\;
{\bf k}\in \mathbb{Z}^3\}.
$$
Define $D^\varepsilon=D\backslash S^\varepsilon$. Then $D^\varepsilon$
is a periodically perforated domain with holes of the same size as
period $\varepsilon$. Notice that the holes are assumed to have no
intersection with the boundary $\partial D$, which implies that
$\partial D^\varepsilon=\partial D\bigcup
\partial S^\varepsilon$. See Fig.\ref{f1}. This assumption is only needed to avoid
technicalities and the results of our paper will remain valid
without the assumption \cite{AMN}.
\par
\par
In the following, we introduce some other notations. Define
$Y^*=Y\backslash \overline{S}$ and $\nu=\frac{|Y^*|}{|Y|}$, with
$|Y|$ and $|Y^*|$ the Lebesgue measure of $Y$ and $Y^*$
respectively. Denote the indicator function $\chi$ as follows
$$
\chi(Y)= \left\{
\begin{array}{ll}
1,&\quad \hbox{on} \quad Y^*,\\
0,&\quad \hbox{on}\quad S,
\end{array}
\right. \quad \hbox{and} \quad \chi(D)= \left\{
\begin{array}{ll}
1,&\quad \hbox{on} \quad D^\varepsilon\\
0,&\quad \hbox{on}\quad S^\varepsilon.
\end{array}
\right.
$$
We also denote $\tilde{u}$ to be the zero extension to the whole
domain $D$ for any function $u$ defined on the domain
$D^\varepsilon$ as follows
$$
\tilde{u}= \left\{
\begin{array}{ll}
u,&\quad \hbox{on} \quad D^\varepsilon,\\
0,&\quad \hbox{on}\quad S^\varepsilon.
\end{array}
\right.
$$
\par
In addition, let the Wiener processes $W_1(t)$ and $W_2(t)$,
defined on a complete probability space $(\Omega, \mathcal{F},
\mathbb{P})$ with the filtration $\{\mathcal{F}_t\}_{t\in
\mathbb{R}}$, be the two-sided in time with values in
$L^2(D)$. Furthermore, assume that $W_1(t)$ and $W_2(t)$ are
 independent and that their covariance operators,   $Q_1$ and
$Q_2$, are   symmetric nonnegative operators  satisfying
$TrQ_1<+\infty$ and $TrQ_2<+\infty$, respectively. Their expansions are given as
follows
$$
\begin{array}{l}
 W_1(t)=\sum\limits_{i=1}^{+\infty}\sqrt{\alpha_{1i}}\beta_{1i} e_i,
\quad \hbox{with}\quad Q_1e_i=\alpha_{1i} e_i,\\
 W_2(t)=\sum\limits_{i=1}^{+\infty}\sqrt{\alpha_{2i}}\beta_{2i} e_i,
\quad \hbox{with}\quad Q_2e_i=\alpha_{2i} e_i,
\end{array}
$$
where $\{e_{i}\}_{i\in \mathbb{N}}$ is an orthonormal bases of
$L^2(D)$,  $\alpha_{1i}$ and $\alpha_{2i}$ are eigenvalues of $Q_1, Q_2$ respectively, and $\{\beta_{1i}\}_{i\in \mathbb{N}}$ and
$\{\beta_{2i}\}_{i\in \mathbb{N}}$ are two sequences of mutually
independent (two-sided in time) standard scalar Wiener process on
the probability space $(\Omega, \mathcal{F}, \mathbb{P})$.

\renewcommand{\theequation}{\thesection.\arabic{equation}}
\setcounter{equation}{0}

\section{Microscopic Model}

\quad\quad Write Equation (\ref{Eq1}) in the It$\hat{o}$ form as
follows

\begin{equation}\label{Eq2}
\left\{
\begin{array}{ll}
du^\varepsilon=v^\varepsilon dt &\quad \hbox{in}\;
D^\varepsilon\times
[0, \tau^*),\\
dv^\varepsilon=(\bigtriangleup
u^\varepsilon-u^\varepsilon-v^\varepsilon+\sin u^\varepsilon)dt+dW_1
&\quad \hbox{in}\; D^\varepsilon\times
[0, \tau^*),\\
d\delta^\varepsilon=\theta^\varepsilon dt & \quad \hbox{on}\; \partial S^\varepsilon\times [0, \tau^*), \\
d\theta^\varepsilon=(-\frac{1}{\varepsilon^2}\theta^\varepsilon
-\delta^\varepsilon-v^\varepsilon)dt+dW_2 & \quad \hbox{on}\;
\partial S^\varepsilon\times [0, \tau^*),\\
u^\varepsilon=v^\varepsilon=0,&\quad \hbox{on}\; \partial D\times[0,
\tau^*),\\
\delta^\varepsilon_t=\frac{\partial u^\varepsilon}{\partial
{\bf n}}& \quad \hbox{on}\;
\partial S^\varepsilon\times [0, \tau^*),\\
\end{array}
\right.
\end{equation}

We supplement Equation (\ref{Eq2}) with the initial data
\begin{equation}\label{Eq3}
u^\varepsilon(0)=u_0, v^\varepsilon(0)=v_0,
\delta^\varepsilon(0)=\delta_0, \theta^\varepsilon(0)=\theta_0,
\end{equation}
which are $\mathcal{F}_0$-measurable.
\par
Now define
$$
A^\varepsilon=\left(
\begin{array}{cccc}
0&\; I&\; 0&\; 0\\
\bigtriangleup-I&\; -I&\; 0&\; 0\\
0&\; 0&\; 0&\; I\\
0&\; -I&\; -I&\; -\frac{1}{\varepsilon^2}I
\end{array}
\right),  F^\varepsilon(U^\varepsilon)=\left(
\begin{array}{l}
0\\
\sin u^\varepsilon\\
0\\
0
\end{array}
\right), W=\left(
\begin{array}{l}
0\\
W_1\\
0\\
W_2
\end{array}
\right).
$$
Let $U^\varepsilon:=(u^\varepsilon, v^\varepsilon,
\delta^\varepsilon, \theta^\varepsilon)^T$ be in the space
$$
\mathcal{H}_\varepsilon: =\{U^\varepsilon\in
H_\varepsilon^1(D^\varepsilon)\times
L_\varepsilon^2(D^\varepsilon)\times L^2(\partial
S^\varepsilon)\times L^2(\partial S^\varepsilon)|\; \frac{\partial
u^\varepsilon}{\partial {\bf n}}=\theta^\varepsilon\; \hbox{on} \;
\partial S^\varepsilon \},$$ with
$$
\|U^\varepsilon\|_{\mathcal{H}_\varepsilon}^2=\|u^\varepsilon\|_{H_\varepsilon^1(D^\varepsilon)}^2
+\|v\|_{L_\varepsilon^2(D^\varepsilon)}^2+\|\delta^\varepsilon\|_{L^2(\partial
S^\varepsilon)}^2+\|\theta\|_{L^2(\partial S^\varepsilon)}^2,
$$
where $H^1_\varepsilon(D^\varepsilon)$ and
$L^2_\varepsilon(D^\varepsilon)$ denote the space
$H^1(D^\varepsilon)$ and $L^2(D^\varepsilon)$ vanishing on $\partial
D$, respectively. The superscript ``$T$" denotes the transpose for
the matrix.
\par
Thus Equation (\ref{Eq2})-(\ref{Eq3}) can be rewritten as
\begin{equation}\label{Eq4}
\left\{
\begin{array}{l}
dU^\varepsilon=A^\varepsilon U^\varepsilon dt+F^\varepsilon(U^\varepsilon)dt+ dW(t),\\
U^\varepsilon(0)=U_0^\varepsilon=(u_0,v_0,\delta_0,\theta_0)^T.
\end{array}
\right.
\end{equation}
\par
For the Cauchy problem (\ref{Eq4}), it follows from Frigeri
\cite{Frigeri} that the operator $A^\varepsilon$ generates a
strongly continuous semigroup $S^\varepsilon(t)=\{e^{A^\varepsilon
t}\}_{t\geq 0}$ on $\mathcal{H}_\varepsilon$. Then the solution of
Equation (\ref{Eq4}) can be written in the mild sense
\begin{equation}\label{E-m}
U^\varepsilon (t)=S^\varepsilon(t)
U^\varepsilon(0)+\int_0^tS^\varepsilon(t-s)F^\varepsilon
(U^\varepsilon(s))ds+\int_0^tS^\varepsilon(t-s) dW(s).
\end{equation}
Furthermore, the variational formulation is
\begin{equation}\label{E-v}
\begin{array}{ll}
&\int_0^{\tau^*}\int_{D^\varepsilon}u^\varepsilon_{tt}\varphi dxdt+
\int_0^{\tau^*}\int_{D^\varepsilon}u^\varepsilon_{t}\varphi dxdt
+\int_0^{\tau^*}\int_{D^\varepsilon}\bigtriangledown u^\varepsilon
\bigtriangledown\varphi dxdt
+\int_0^{\tau^*}\int_{D^\varepsilon}u^\varepsilon\varphi dxdt\\
&-\int_0^{\tau^*}\int_{D^\varepsilon}\sin u^\varepsilon\varphi dxdt
+\varepsilon^2\int_0^{\tau^*}\int_{\partial
S^\varepsilon}\delta^\varepsilon_{tt}\varphi dxdt
+\varepsilon^2\int_0^{\tau^*}\int_{\partial
S^\varepsilon}\delta^\varepsilon\varphi dxdt\\
=& \int_0^{\tau^*}\int_{D^\varepsilon}\dot{W}_1\varphi dxdt -
\varepsilon^2\int_0^{\tau^*}\int_{\partial
S^\varepsilon}u^\varepsilon_t\varphi dxdt
+\varepsilon^2\int_0^{\tau^*}\int_{\partial
S^\varepsilon}\dot{W}_2\varphi dxdt,
\end{array}
\end{equation}
for any $\varphi\in C_0^\infty([0, \tau^*)\times D^\varepsilon)$.
\par
{\bf Proposition 3.1} {(\bf Local well-posedness)}\quad {\it Let the
initial datum $U_0^\varepsilon$ be a $ \mathcal{F}_0$-measurable
random variable with value in $\mathcal{H}_\varepsilon$. Then the
Cauchy problem (\ref{Eq4}) has a unique local mild solution
$U^\varepsilon(t)$ in $C([0, \tau^*), \mathcal{H}_\varepsilon)$,
where $\tau^*$ is a stopping time depending on $U_0^\varepsilon$ and
$\omega$. Moreover, the mild solution $U^\varepsilon(t)$ is also a
weak solution in the following sense
\begin{equation}
\langle U^\varepsilon (t), \phi\rangle_{\mathcal{H}_\varepsilon}=
\langle U^\varepsilon(0), \phi \rangle_{\mathcal{H}_\varepsilon}
+\int_0^t \langle A^\varepsilon U^\varepsilon (s),
\phi\rangle_{\mathcal{H}_\varepsilon}ds+ \int_0^t \langle
F^\varepsilon(U^\varepsilon (s)),
\phi\rangle_{\mathcal{H}_\varepsilon}ds+\int_0^t \langle dW(s),\phi
\rangle_{\mathcal{H}_\varepsilon}
\end{equation}
for any $t\in[0, \tau^*)$ and $\phi\in \mathcal{H}_\varepsilon$. }
\par
{\bf Proof.}\quad We first define a cut-off function as follows. For
any positive parameter $R$, let $\eta_R(\cdot)$ be a positive real
valued $C^\infty$-function on $[0, +\infty)$ such that
$$
\eta_R(s)=\left\{
\begin{array}{ll}
1, & \quad \hbox{for}\quad 0\leq s\leq \frac{R}{2},\\
\in(0, 1), &\quad \hbox{for}\quad \frac{R}{2}< s\leq R,\\
0, &\quad \hbox{for}\quad R< s< +\infty.
\end{array}
\right.
$$
\par
Then the truncated system of the Cauchy problem (\ref{Eq4}) is
defined as follows
\begin{equation}\label{Eq5}
\left\{
\begin{array}{l}
dU^\varepsilon=A^\varepsilon U^\varepsilon dt+F^\varepsilon_R(U^\varepsilon)dt+ dW(t),\\
U^\varepsilon(0)=U_0,
\end{array}
\right.
\end{equation}
where $F^\varepsilon_R(U^\varepsilon)=(0,
\eta_R(\|U^\varepsilon\|_{\mathcal{H}_\varepsilon}^2)\sin
u^\varepsilon, 0, 0)^{T}$.
\par
In the meantime, we easily examine that
$F^\varepsilon_R(U^\varepsilon)$ satisfies the sublinear growth and
the Lipschitz continuity as in Chen and Zhang \cite{Chen2}.
Therefore, according to Theorem 7.4 of Da Prato and Zabczyk
\cite{DZ1}, the truncated system (\ref{Eq5}) has a unique mild
solution $U^\varepsilon_R(t)$ in $\mathcal{H}_\varepsilon$ for each
fixed positive $R$.
\par
Define a stopping time
\begin{equation}\label{Eq7}
\tau_R:=\inf\{t:
\|U^\varepsilon\|_{\mathcal{H}_\varepsilon}^2>\frac{R}{2}\}.
\end{equation}
We have $U^\varepsilon(t)=U^\varepsilon_R(t)$ as $t<\tau_R$. Also
from Da Prato and Zabczyk \cite{DZ1}, the path   $t\to
U^\varepsilon(t)$ is continuous. Let
$\tau^*=\lim\limits_{R\to+\infty}\tau_R$. Then $U^\varepsilon(t)$ is
the unique local solution of the Cauchy problem (\ref{Eq4}) with
lifespan $\tau^*$. Furthermore, applying the stochastic Fubini
theorem, it can be verified that the local mild solution is also the
weak solution. The proof is complete. \hfill$\blacksquare$
\par

Because the energy relation of this stochastic system does not
directly imply the a priori estimate of the solutions, we will
introduce a pseudo energy argument (see Chow\cite{Chow2} and Chen
and Zhang \cite{Chen2}) to establish the a priori estimate of the
solutions for the Cauchy problem (\ref{Eq4}). Furthermore, applying
the a priori estimate, we could obtain the global existence and
almost sure boundedness of solutions, which further implies the
tightness of distribution of solutions.
\par
For   a real parameter $r$ in $(0,1)$, we define
\begin{equation}\label{T-Psolution}
v_r^\varepsilon=v^\varepsilon+ru^\varepsilon\quad \hbox{and}\quad
\theta_r^\varepsilon=\theta^\varepsilon+r\delta^\varepsilon,
\end{equation}
with $(u^\varepsilon, v^\varepsilon, \delta^\varepsilon,
\theta^\varepsilon)^T$ being the solution of the Cauchy problem
(\ref{Eq2})-(\ref{Eq3}). Then the solution $U_r^\varepsilon=
(u^\varepsilon, v_r^\varepsilon, \delta^\varepsilon,
\theta_r^\varepsilon)^T \in \mathcal{H}_\varepsilon$ satisfies the
following equation
\begin{equation}\label{Eq-P}
\left\{
\begin{array}{ll}
du^\varepsilon=(v_r^\varepsilon-r u^\varepsilon)dt & \hbox{in}\;
D^\varepsilon\times
[0, \tau^*),\\
dv_r^\varepsilon=(\bigtriangleup
u^\varepsilon-(1-r+r^2)u^\varepsilon-(1-r)v_r^\varepsilon+\sin
u^\varepsilon)dt+dW_1 & \hbox{in}\; D^\varepsilon\times
[0, \tau^*),\\
d\delta^\varepsilon=(\theta_r^\varepsilon-r\delta^\varepsilon)dt &  \hbox{on}\;
\partial S^\varepsilon \times [0, \tau^*), \\
d\theta_r^\varepsilon=(-(\frac{1}{\varepsilon^2}-r)\theta_r^\varepsilon
-(1-\frac{r}{\varepsilon^2}+r^2)\delta^\varepsilon-v_r^\varepsilon+r
u^\varepsilon)dt+dW_2 &  \hbox{on}\; \partial S^\varepsilon\times [0, \tau^*),\\
u^\varepsilon=v_r^\varepsilon=0& \hbox{on}\; \partial D\times[0,
\tau^*),\\
\delta^\varepsilon_t=\frac{\partial u^\varepsilon}{\partial {\bf
n}}& \hbox{on}\;
\partial S^\varepsilon \times [0, \tau^*),\\
u^\varepsilon(0)=u_0, v_r^\varepsilon(0)=v_0+r u_0:=v_{r0}& \hbox{in}\;  D^\varepsilon,\\
\delta^\varepsilon(0)=\delta_0,
\theta_r^\varepsilon(0)=\theta_0+r\delta_0:=\theta_{r0} &
\hbox{on}\;
\partial S^\varepsilon.
\end{array}
\right.
\end{equation}
\par

Define the pseudo energy functional $\mathcal{E}_r^\varepsilon(t)$
of the Cauchy problem (\ref{Eq4}) as follows
$$
\begin{array}{lll}
\mathcal{E}_r^\varepsilon(t)&:=&\|v_r^\varepsilon(t)\|_{L_\varepsilon^2(D^\varepsilon)}^2+\|\bigtriangledown
u^\varepsilon(t)\|_{L_\varepsilon^2(D^\varepsilon)}^2+(1-r+r^2)\|u^\varepsilon(t)\|_{L_\varepsilon^2(D^\varepsilon)}^2\\
&& +\|\theta_r^\varepsilon(t)\|_{L^2(\partial S^\varepsilon)}^2
+(1-\frac{r}{\varepsilon^2}+r^2)\|\delta^\varepsilon(t)\|_{L^2(\partial
S^\varepsilon)}^2
+4\|\cos \frac{u^\varepsilon(t)}{2}\|_{L_\varepsilon^2(D^\varepsilon)}^2\\
&& +2r \langle
u^\varepsilon(t),\delta^\varepsilon(t)\rangle_{L^2(\partial
S^\varepsilon)}.
\end{array}
$$
\par

{\bf Proposition 3.2}\quad{\it Let the initial data
$U_r^\varepsilon(0)$ be a $ \mathcal{F}_0$-measurable random
variable in $L^2(\Omega, \mathcal{H}_\varepsilon)$. Then for any
time $t\in [0, \tau^*)$, we have
\begin{equation}\label{Eq-3.2-01}
\begin{array}{ll}
\mathcal{E}_r^\varepsilon(t)=&
\mathcal{E}_r^\varepsilon(0)-\int_0^t[2(1-r)\|v_r^\varepsilon\|_{L_\varepsilon^2(D^\varepsilon)}^2+2r\|\bigtriangledown
u^\varepsilon\|_{L_\varepsilon^2(D^\varepsilon)}^2+2r(1-r+r^2)\|u^\varepsilon\|_{L_\varepsilon^2(D^\varepsilon)}^2\\
&\quad\quad\quad\quad
+2(\frac{1}{\varepsilon^2}-r)\|\theta_r^\varepsilon\|_{L^2(\partial
S^\varepsilon)}^2+2(1-\frac{r}{\varepsilon^2}
+r^2)r\|\delta^\varepsilon\|_{L^2(\partial
S^\varepsilon)}^2]ds\\
&+2r\int_0^t\langle u^\varepsilon,\sin
u^\varepsilon\rangle_{L_\varepsilon^2(D^\varepsilon)}ds+4r\int_0^t
\langle u^\varepsilon,\theta_r^\varepsilon\rangle_{L^2(\partial
S^\varepsilon)}ds-4r^2\int_0^t
\langle u^\varepsilon,\delta^\varepsilon\rangle_{L^2(\partial S^\varepsilon)}ds\\
&+\int_0^t\langle 2v_r^\varepsilon,
dW_1(s)\rangle_{L_\varepsilon^2(D^\varepsilon)}+\int_0^t\langle
2\theta_r^\varepsilon,
dW_2(s)\rangle_{L^2(\partial S^\varepsilon)}\\
&+tTrQ_1+tTrQ_2.
\end{array}
\end{equation}
Moreover,
\begin{equation}\label{Eq-3.2-02}
\begin{array}{ll}
\mathbb{E} \mathcal{E}_r^\varepsilon(t)=&
\mathbb{E}\mathcal{E}_r^\varepsilon(0)\\
&-\int_0^t[2(1-r)\mathbb{E}\|v_r^\varepsilon\|_{L_\varepsilon^2(D^\varepsilon)}^2+2r\mathbb{E}\|\bigtriangledown
u^\varepsilon\|_{L_\varepsilon^2(D^\varepsilon)}^2
+2r(1-r+r^2)\mathbb{E}\|u^\varepsilon\|_{L_\varepsilon^2(D^\varepsilon)}^2\\
&\quad\quad+2(\frac{1}{\varepsilon^2}-r)\mathbb{E}\|\theta_r^\varepsilon\|_{L^2(\partial
S^\varepsilon)}^2+2(1-\frac{r}{\varepsilon^2}
+r^2)r\mathbb{E}\|\delta^\varepsilon\|_{L^2(\partial
S^\varepsilon)}^2]ds\\
&+2r\int_0^t\mathbb{E}\langle u^\varepsilon,\sin
u^\varepsilon\rangle_{L_\varepsilon^2(D^\varepsilon)}ds+4r\int_0^t
\mathbb{E}\langle
u^\varepsilon,\theta_r^\varepsilon\rangle_{L^2(\partial
S^\varepsilon)}ds\\
&-4r^2\int_0^t \mathbb{E}\langle
u^\varepsilon,\delta^\varepsilon\rangle_{L^2(\partial
S^\varepsilon)}ds+tTrQ_1+tTrQ_2.
\end{array}
\end{equation}
}

\par
{\bf  Proof.}\quad  First, we examine the second equation of
(\ref{Eq-P}). Put $M(v_r^\varepsilon):=\int_{D^\varepsilon}
|v_r^\varepsilon|^2dx$. Then from It$\hat{o}$ formula, we
deduce that
\begin{equation}\label{Eq-3.2-1}
\begin{array}{ll}
M(v_r^\varepsilon(t))=& M(v_r^\varepsilon(0))+\int_0^t\langle
M^{\prime}(v_r^\varepsilon),
dW_1(s)\rangle_{L_\varepsilon^2(D^\varepsilon)}
+\int_0^t\frac{1}{2}Tr[M^{\prime\prime}(v_r^\varepsilon)Q_1^{\frac{1}{2}}(Q_1^{\frac{1}{2}})^*]ds\\
&+ \int_0^t\langle M^{\prime}(v_r^\varepsilon),(\bigtriangleup
u^\varepsilon-(1-r+r^2)u^\varepsilon-(1-r)v_r^\varepsilon+\sin
u^\varepsilon)\rangle_{L_\varepsilon^2(D^\varepsilon)}ds,
\end{array}
\end{equation}
with   $M^\prime(v_r^\varepsilon)=2v_r^\varepsilon$ and
$M^{\prime\prime}(v_r^\varepsilon)=2\varphi$ for any $\varphi$ in
$L_\varepsilon^2(D^\varepsilon)$. After some calculations, we   get
that
\begin{equation}\label{Eq-3.2-2}
\begin{array}{ll}
&\langle M^{\prime}(v_r^\varepsilon),(\bigtriangleup
u^\varepsilon-(1-r+r^2)u^\varepsilon-(1-r)v_r^\varepsilon+\sin
u^\varepsilon)\rangle_{L_\varepsilon^2(D^\varepsilon)}\\
=& -\frac{d}{ds}[\|\bigtriangledown
u^\varepsilon\|_{L_\varepsilon^2(D^\varepsilon)}^2+(1-r+r^2)\|u^\varepsilon\|_{L_\varepsilon^2(D^\varepsilon)}^2+4\|\cos
\frac{u^\varepsilon}{2}\|_{L_\varepsilon^2(D^\varepsilon)}^2]\\
&-[2r\|\bigtriangledown
u^\varepsilon\|_{L_\varepsilon^2(D^\varepsilon)}^2+2r(1-r+r^2)\|u^\varepsilon\|_{L_\varepsilon^2(D^\varepsilon)}^2
+2(1-r)\|v_r^\varepsilon\|_{L_\varepsilon^2(D^\varepsilon)}^2]\\
&+2\langle v_r^\varepsilon,\frac{\partial u^\varepsilon}{\partial
{\bf n}}\rangle_{L^2(\partial S^\varepsilon)}+2r\langle
u^\varepsilon,\sin
u^\varepsilon\rangle_{L_\varepsilon^2(D^\varepsilon)}.
\end{array}
\end{equation}
It immediately follows from (\ref{Eq-3.2-1}) and (\ref{Eq-3.2-2})
that
\begin{equation}\label{Eq-3.2-3}
\begin{array}{ll}
&\|v_r^\varepsilon(t)\|_{L_\varepsilon^2(D^\varepsilon)}^2+\|\bigtriangledown
u^\varepsilon(t)\|_{L_\varepsilon^2(D^\varepsilon)}^2+(1-r+r^2)
\|u^\varepsilon(t)\|_{L_\varepsilon^2(D^\varepsilon)}^2+4\|\cos
\frac{u^\varepsilon(t)}{2}\|_{L_\varepsilon^2(D^\varepsilon)}^2\\
=&
\|v_r^\varepsilon(0)\|_{L_\varepsilon^2(D^\varepsilon)}^2+\|\bigtriangledown
u^\varepsilon(0)\|_{L_\varepsilon^2(D^\varepsilon)}^2+(1-r+r^2)\|u^\varepsilon(0)\|_{L_\varepsilon^2(D^\varepsilon)}^2+4\|\cos
\frac{u^\varepsilon(0)}{2}\|_{L_\varepsilon^2(D^\varepsilon)}^2\\
&-\int_0^t[2(1-r)\|v_r^\varepsilon\|_{L_\varepsilon^2(D^\varepsilon)}^2+2r\|\bigtriangledown
u^\varepsilon\|_{L_\varepsilon^2(D^\varepsilon)}^2+2r(1-r+r^2)\|u^\varepsilon\|_{L_\varepsilon^2(D^\varepsilon)}^2]ds\\
&+2\int_0^t\langle v_r^\varepsilon,\frac{\partial
u^\varepsilon}{\partial {\bf n}}\rangle_{L^2(\partial
S^\varepsilon)}ds+2r\int_0^t\langle u^\varepsilon,\sin
u^\varepsilon\rangle_{L_\varepsilon^2(D^\varepsilon)}ds\\
&+\int_0^t\langle 2v_r^\varepsilon,
dW_1(s)\rangle_{L_\varepsilon^2(D^\varepsilon)}+tTrQ_1.
\end{array}
\end{equation}
\par
Second, we examine the fourth equation of (\ref{Eq-P}) and
$M(\theta_r^\varepsilon)=\int_{\partial S^\varepsilon}
|\theta_r^\varepsilon|^2dx$. Note that
\begin{equation}\label{Eq-3.2-4}
\begin{array}{lll}
M(\theta_r^\varepsilon(t))&=&
M(\theta_r^\varepsilon(0))+\int_0^t\langle
M^{\prime}(\theta_r^\varepsilon), dW_2(s)\rangle_{L^2(\partial
S^\varepsilon)}
+\int_0^t\frac{1}{2}Tr[M^{\prime\prime}(\theta_r^\varepsilon)Q_2^{\frac{1}{2}}(Q_2^{\frac{1}{2}})^*]ds\\
&&+ \int_0^t\langle
M^{\prime}(\theta_r^\varepsilon),(-(\frac{1}{\varepsilon^2}-r)\theta_r^\varepsilon-
(1-\frac{r}{\varepsilon^2}+r^2)\delta^\varepsilon-v_r^\varepsilon+r
u^\varepsilon)\rangle_{L^2(\partial S^\varepsilon)}ds,
\end{array}
\end{equation}
with   $M^\prime(\theta_r^\varepsilon)=2\theta_r^\varepsilon$ and
$M^{\prime\prime}(\theta_r^\varepsilon)=2\phi$ for any $\phi$ in
$L^2(\partial S^\varepsilon)$. After some calculations, we   conclude
that
\begin{equation}\label{Eq-3.2-5}
\begin{array}{ll}
&\langle
M^{\prime}(\theta_r^\varepsilon),(-(\frac{1}{\varepsilon^2}-r)\theta_r^\varepsilon-
(1-\frac{r}{\varepsilon^2}+r^2)\delta^\varepsilon-v_r^\varepsilon+r
u^\varepsilon)\rangle_{L^2(\partial S^\varepsilon)}\\
=&-(1-\frac{r}{\varepsilon^2}+r^2)\frac{d}{ds}\|\delta^\varepsilon\|_{L^2(\partial
S^\varepsilon)}^2
-2(1-\frac{r}{\varepsilon^2}+r^2)r\|\delta^\varepsilon\|_{L^2(\partial
S^\varepsilon)}^2-2(\frac{1}{\varepsilon^2}-r)\|\theta_r^\varepsilon\|_{L^2(\partial
S^\varepsilon)}^2\\
&-2\langle \frac{\partial u^\varepsilon}{\partial {\bf
n}},v_r^\varepsilon\rangle_{L^2(\partial S^\varepsilon)}-2r\langle
\delta^\varepsilon,v_r^\varepsilon\rangle_{L^2(\partial
S^\varepsilon)} +2r\langle \theta_r^\varepsilon,
u^\varepsilon\rangle_{L^2(\partial S^\varepsilon)}.
\end{array}
\end{equation}
It   follows from (\ref{Eq-3.2-4}) and (\ref{Eq-3.2-5})
that
\begin{equation}\label{Eq-3.2-6}
\begin{array}{ll}
&\|\theta_r^\varepsilon(t)\|_{L^2(\partial S^\varepsilon)}^2
+(1-\frac{r}{\varepsilon^2}+r^2)\|\delta^\varepsilon(t)\|_{L^2(\partial S^\varepsilon)}^2\\
=&\|\theta^\varepsilon(0)\|_{L^2(\partial S^\varepsilon)}^2
+(1-\frac{r}{\varepsilon^2}+r^2)\|\delta^\varepsilon(0)\|_{L^2(\partial S^\varepsilon)}^2\\
&-\int_0^t[2(1-\frac{r}{\varepsilon^2}+r^2)r\|\delta^\varepsilon\|_{L^2(\partial
S^\varepsilon)}^2+2(\frac{1}{\varepsilon^2}-r)\|\theta_r^\varepsilon\|_{L^2(\partial
S^\varepsilon)}^2]ds\\
&-2\int_0^t\langle \frac{\partial u^\varepsilon}{\partial {\bf
n}},v_r^\varepsilon\rangle_{L^2(\partial
S^\varepsilon)}ds-2r\int_0^t\langle
\delta^\varepsilon,v_r^\varepsilon\rangle_{L^2(\partial
S^\varepsilon)}ds +2r\int_0^t\langle \theta_r^\varepsilon,
u^\varepsilon\rangle_{L^2(\partial S^\varepsilon)}ds\\
&+\int_0^t\langle 2\theta_r^\varepsilon,
dW_2(s)\rangle_{L^2(\partial S^\varepsilon)}+tTrQ_2.
\end{array}
\end{equation}
\par
Thus, from (\ref{Eq-3.2-3}) and (\ref{Eq-3.2-6}), we have
\begin{equation}\label{Eq-3.2-7}
\begin{array}{ll}
&\|v_r^\varepsilon(t)\|_{L_\varepsilon^2(D^\varepsilon)}^2+\|\bigtriangledown
u^\varepsilon(t)\|_{L_\varepsilon^2(D^\varepsilon)}^2+(1-r+r^2)
\|u^\varepsilon(t)\|_{L_\varepsilon^2(D^\varepsilon)}^2+\|\theta_r^\varepsilon(t)\|_{L^2(\partial S^\varepsilon)}^2\\
&+(1-\frac{r}{\varepsilon^2}+r^2)\|\delta^\varepsilon(t)\|_{L^2(\partial
S^\varepsilon)}^2+4\|\cos
\frac{u^\varepsilon(t)}{2}\|_{L_\varepsilon^2(D^\varepsilon)}^2\\
=&
\|v_r^\varepsilon(0)\|_{L_\varepsilon^2(D^\varepsilon)}^2+\|\bigtriangledown
u^\varepsilon(0)\|_{L_\varepsilon^2(D^\varepsilon)}^2+(1-r+r^2)\|u^\varepsilon(0)\|_{L_\varepsilon^2(D^\varepsilon)}^2
+\|\theta^\varepsilon(0)\|_{L^2(\partial
S^\varepsilon)}^2\\
&+(1-\frac{r}{\varepsilon^2}+r^2)\|\delta^\varepsilon(0)\|_{L^2(\partial
S^\varepsilon)}^2+4\|\cos
\frac{u^\varepsilon(0)}{2}\|_{L_\varepsilon^2(D^\varepsilon)}^2\\
&-\int_0^t[2(1-r)\|v_r^\varepsilon\|_{L_\varepsilon^2(D^\varepsilon)}^2+2r\|\bigtriangledown
u^\varepsilon\|_{L_\varepsilon^2(D^\varepsilon)}^2+2r(1-r+r^2)\|u^\varepsilon\|_{L_\varepsilon^2(D^\varepsilon)}^2\\
&+2(\frac{1}{\varepsilon^2}-r)\|\theta_r^\varepsilon\|_{L^2(\partial
S^\varepsilon)}^2+2(1-\frac{r}{\varepsilon^2}+r^2)r\|\delta^\varepsilon\|_{L^2(\partial
S^\varepsilon)}^2
]ds\\
&+2r\int_0^t\langle u^\varepsilon,\sin
u^\varepsilon\rangle_{L_\varepsilon^2(D^\varepsilon)}ds-2r\int_0^t\langle
\delta^\varepsilon,v_r^\varepsilon\rangle_{L^2(\partial
S^\varepsilon)}ds +2r\int_0^t\langle \theta_r^\varepsilon,
u^\varepsilon\rangle_{L^2(\partial S^\varepsilon)}ds\\
&+\int_0^t\langle 2v_r^\varepsilon,
dW_1(s)\rangle_{L_\varepsilon^2(D^\varepsilon)}+\int_0^t\langle
2\theta_r^\varepsilon, dW_2(s)\rangle_{L^2(\partial
S^\varepsilon)}+tTrQ_1+tTrQ_2.
\end{array}
\end{equation}
Meanwhile, we note that
$$
\begin{array}{ll}
&2r \langle u^\varepsilon(t),\delta^\varepsilon(t)\rangle_{L^2(\partial S^\varepsilon)}\\
= &2r \langle
u^\varepsilon(0),\delta^\varepsilon(0)\rangle_{L^2(\partial
S^\varepsilon)}+2r\int_0^t \langle
(u^\varepsilon)_s,\delta^\varepsilon \rangle_{L^2(\partial
S^\varepsilon)}ds+2r\int_0^t \langle
u^\varepsilon,(\delta^\varepsilon)_s \rangle_{L^2(\partial S^\varepsilon)}ds\\
= &2r \langle
u^\varepsilon(0),\delta^\varepsilon(0)\rangle_{L^2(\partial
S^\varepsilon)}+2r\int_0^t \langle v^\varepsilon,\delta^\varepsilon
\rangle_{L^2(\partial S^\varepsilon)}ds+2r\int_0^t \langle
u^\varepsilon,\theta^\varepsilon \rangle_{L^2(\partial S^\varepsilon)}ds\\
= &2r \langle
u^\varepsilon(0),\delta^\varepsilon(0)\rangle_{L^2(\partial
S^\varepsilon)}+2r\int_0^t \langle v_r^\varepsilon-r
u^\varepsilon,\delta^\varepsilon \rangle_{L^2(\partial
S^\varepsilon)}ds+2r\int_0^t \langle
u^\varepsilon,\theta_r^\varepsilon-r\delta^\varepsilon \rangle_{L^2(\partial S^\varepsilon)}ds\\
= &2r \langle
u^\varepsilon(0),\delta^\varepsilon(0)\rangle_{L^2(\partial
S^\varepsilon)}+2r\int_0^t \langle
v_r^\varepsilon,\delta^\varepsilon \rangle_{L^2(\partial
S^\varepsilon)}ds+2r\int_0^t \langle
u^\varepsilon,\theta_r^\varepsilon\rangle_{L^2(\partial S^\varepsilon)}ds\\
&-4r^2\int_0^t \langle
u^\varepsilon,\delta^\varepsilon\rangle_{L^2(\partial
S^\varepsilon)}ds,
\end{array}
$$
which implies that
\begin{equation}\label{Eq-3.2-8}
\begin{array}{ll}
&-2r\int_0^t \langle v_r^\varepsilon,\delta^\varepsilon
\rangle_{L^2(\partial S^\varepsilon)}ds
\\
=&-2r \langle
u^\varepsilon(t),\delta^\varepsilon(t)\rangle_{L^2(\partial
S^\varepsilon)}+2r \langle
u^\varepsilon(0),\delta^\varepsilon(0)\rangle_{L^2(\partial
S^\varepsilon)}+2r\int_0^t \langle
u^\varepsilon,\theta_r^\varepsilon\rangle_{L^2(\partial S^\varepsilon)}ds\\
&-4r^2\int_0^t \langle
u^\varepsilon,\delta^\varepsilon\rangle_{L^2(\partial
S^\varepsilon)}ds.
\end{array}
\end{equation}
Then it   follows from (\ref{Eq-3.2-7}) and
(\ref{Eq-3.2-8}) that (\ref{Eq-3.2-01}) and (\ref{Eq-3.2-02}) hold.
\hfill$\blacksquare$

\par

{\bf Proposition 3.3}\quad {\it Let the initial datum
$U_r^\varepsilon(0)$ be a $ \mathcal{F}_0$-measurable random
variable in $L^2(\Omega, \mathcal{H}_\varepsilon)$. Then for any
time $t\in [0, \tau^*)$, and a sufficient small $r$ in $(0, 1)$,
there exists a positive constant $C$ such that
\begin{equation}\label{Eq-3.3-0}
\mathbb{E}\|U_r^\varepsilon(t)\|_{\mathcal{H}_\varepsilon}^2 \leq
C\mathbb{E}
\mathcal{E}_r^\varepsilon(0)-C\int_0^t[\mathbb{E}\|U_r^\varepsilon(s)\|_{\mathcal{H}_\varepsilon}^2]ds+C[tTrQ_1
+tTrQ_2+t].
\end{equation}
}
\par
{\bf Proof.}\quad On the one hand, it follows from the Cauchy
inequality and the trace inequality that there exists a positive
constant $C_{TI}>0$ (here and hereafter $C_{TI}$ denotes the
positive constant in the trace inequality) such that
$$
\begin{array}{ll}
0 &\leq r\mathbb{E}\|u^\varepsilon(t)\|_{L^2(\partial
S^\varepsilon)}^2+2r \mathbb{E}\langle
u^\varepsilon(t),\delta^\varepsilon(t)\rangle_{L^2(\partial
S^\varepsilon)}+r
\mathbb{E}\|\delta^\varepsilon(t)\|_{L^2(\partial S^\varepsilon)}^2\\
&\leq r
C_{TI}^2\mathbb{E}\|u^\varepsilon(t)\|_{H_\varepsilon^1(D^\varepsilon)}^2+2r
\mathbb{E}\langle
u^\varepsilon(t),\delta^\varepsilon(t)\rangle_{L^2(\partial
S^\varepsilon)}+r \mathbb{E}\|\delta^\varepsilon(t)\|_{L^2(\partial
S^\varepsilon)}^2,
\end{array}
$$
which implies that
\begin{equation}\label{Eq-3.3-1}
\begin{array}{ll}
\mathbb{E} \mathcal{E}_r^\varepsilon(t)\geq
&\mathbb{E}\|v_r^\varepsilon(t)\|_{L_\varepsilon^2(D^\varepsilon)}^2+(1-r
C_{TI}^2)\mathbb{E}\|\bigtriangledown
u^\varepsilon(t)\|_{L_\varepsilon^2(D^\varepsilon)}^2\\
&+(1-r-r
C_{TI}^2+r^2)\mathbb{E}\|u^\varepsilon(t)\|_{L_\varepsilon^2(D^\varepsilon)}^2
+\mathbb{E}\|\theta_r^\varepsilon(t)\|_{L^2(\partial S^\varepsilon)}^2\\
&+(1-\frac{r}{\varepsilon^2}-r+r^2)\mathbb{E}\|\delta(t)\|_{L^2(\partial
S^\varepsilon)}^2.
\end{array}
\end{equation}
\par

On the other hand, it follows from the H$\ddot{o}$lder inequality,
the Young inequality and the trace inequality that
$$
\begin{array}{ll}
\mathbb{E}\langle
u^\varepsilon,\theta_r^\varepsilon\rangle_{L^2(\partial
S^\varepsilon)}&\leq
\mathbb{E}\|u^\varepsilon\|_{L^2(\partial S^\varepsilon)}\cdot\mathbb{E}\|\theta_r^\varepsilon\|_{L^2(\partial S^\varepsilon)}\\
&\leq r
\mathbb{E}\|u^\varepsilon\|_{L^2(\partial S^\varepsilon)}^2+\frac{1}{4r}\mathbb{E}\|\theta_r^\varepsilon\|_{L^2(\partial S^\varepsilon)}^2\\
&\leq r C_{TI}^2
\mathbb{E}\|u^\varepsilon\|_{H_\varepsilon^1(D^\varepsilon)}^2+\frac{1}{4r}\mathbb{E}\|\theta_r^\varepsilon\|_{L^2(\partial
S^\varepsilon)}^2,
\end{array}
$$
which implies that
\begin{equation}\label{Eq-3.3-2}
4r \mathbb{E}\langle
u^\varepsilon,\theta_r^\varepsilon\rangle_{L^2(\partial
S^\varepsilon)}\leq 4r^2 C_{TI}^2\mathbb{E}\|\bigtriangledown
u^\varepsilon\|_{L_\varepsilon^2(D^\varepsilon)}^2+4r^2
C_{TI}^2\mathbb{E}\|u^\varepsilon\|_{L_\varepsilon^2(D^\varepsilon)}^2
+\mathbb{E}\|\theta_r^\varepsilon\|_{L^2(\partial S^\varepsilon)}^2.
\end{equation}
\par
At the same time, it follows from the Cauchy inequality and the
trace inequality that
\begin{equation}\label{Eq-3.3-3}
\begin{array}{ll}
-4r^2 \mathbb{E}\langle
u^\varepsilon,\delta^\varepsilon\rangle_{L^2(\partial
S^\varepsilon)}&\leq 2r^2
\mathbb{E}\|u^\varepsilon\|_{L^2(\partial S^\varepsilon)}^2+2r^2 \mathbb{E}\|\delta^\varepsilon\|_{L^2(\partial S^\varepsilon)}^2\\
&\leq 2r^2C_{TI}^2 \mathbb{E}\|\bigtriangledown
u^\varepsilon\|_{L_\varepsilon^2(D^\varepsilon)}^2+2r^2C_{TI}^2
\mathbb{E}\|u^\varepsilon\|_{L_\varepsilon^2(D^\varepsilon)}^2+2r^2
\mathbb{E}\|\delta^\varepsilon\|_{L^2(\partial S^\varepsilon)}^2.
\end{array}
\end{equation}
Also it follows from the Cauchy inequality that
\begin{equation}\label{Eq-3.3-4}
\begin{array}{ll}
2r\mathbb{E}\langle u^\varepsilon,\sin
u^\varepsilon\rangle_{L_\varepsilon^2(D^\varepsilon)}&\leq r
\mathbb{E}\|u^\varepsilon\|_{L_\varepsilon^2(D^\varepsilon)}^2+ r
\mathbb{E}\|\sin
u^\varepsilon\|_{L_\varepsilon^2(D^\varepsilon)}^2\\
&\leq r
\mathbb{E}\|u^\varepsilon\|_{L_\varepsilon^2(D^\varepsilon)}^2+ C.
\end{array}
\end{equation}
\par
Notice that $\varepsilon \in (0, 1)$. Then it follows from
Proposition 3.2 and (\ref{Eq-3.3-2})-(\ref{Eq-3.3-4}) that
\begin{equation}\label{Eq-3.3-5}
\begin{array}{ll}
\mathbb{E} \mathcal{E}_r^\varepsilon(t)\leq& \mathbb{E}
\mathcal{E}_r^\varepsilon(0)-\int_0^t[2(1-r)\mathbb{E}\|v^\varepsilon\|_{L_\varepsilon^2(D^\varepsilon)}^2
+2r(1-3r C_{TI}^2)\mathbb{E}\|\bigtriangledown
u^\varepsilon\|_{L_\varepsilon^2(D^\varepsilon)}^2\\
&+r[1-2r-6r C_{TI}^2+2r^2]\mathbb{E}\|u^\varepsilon\|_{L_\varepsilon^2(D^\varepsilon)}^2\\
&+(1-2r)\mathbb{E}\|\theta^\varepsilon\|_{L^2(\partial
S^\varepsilon)}^2
+2r(1-r-\frac{r}{\varepsilon^2}+r^2)\mathbb{E}\|\delta^\varepsilon\|_{L^2(\partial S^\varepsilon)}^2]ds\\
&+tTrQ_1+tTrQ_2+Ct.
\end{array}
\end{equation}

Let $r$ be sufficient small in $(0, 1)$ such that
\begin{equation}\label{Eq-3.3-6}
\min\{1-2r, 1-3r C_{TI}^2, 1-r-\frac{r}{\varepsilon^2}+r^2, 1-2r-6r
C_{TI}^2+2r^2, 1-r-r C_{TI}^2+r^2
 \}>0.
\end{equation}

Therefore, from (\ref{Eq-3.3-1}), (\ref{Eq-3.3-5}) and
(\ref{Eq-3.3-6}), there exists a positive constant $C$ such that
$$
\begin{array}{ll}
&\mathbb{E}\|v^\varepsilon(t)\|_{L_\varepsilon^2(D^\varepsilon)}^2+\mathbb{E}\|\bigtriangledown
u^\varepsilon(t)\|_{L_\varepsilon^2(D^\varepsilon)}^2+\mathbb{E}\|u^\varepsilon(t)\|_{L_\varepsilon^2(D^\varepsilon)}^2\\
&+\mathbb{E}\|\theta_r^\varepsilon(t)\|_{L^2(\partial
S^\varepsilon)}^2
+\mathbb{E}\|\delta^\varepsilon(t)\|_{L^2(\partial S^\varepsilon)}^2\\
\leq & C
\mathcal{E}_r^\varepsilon(0)-C\int_0^t[\mathbb{E}\|v^\varepsilon(s)\|_{L_\varepsilon^2(D^\varepsilon)}^2+\mathbb{E}\|\bigtriangledown
u^\varepsilon(s)\|_{L_\varepsilon^2(D^\varepsilon)}^2\\
&+\mathbb{E}\|u^\varepsilon(s)\|_{L_\varepsilon^2(D^\varepsilon)}^2+\mathbb{E}\|\theta_r^\varepsilon(s)\|_{L^2(\partial
S^\varepsilon)}^2
\\
&+\mathbb{E}\|\delta^\varepsilon(s)\|_{L^2(\partial
S^\varepsilon)}^2]ds+C[tTrQ_1 +tTrQ_2+t],
\end{array}
$$
which implies (\ref{Eq-3.3-0}). \hfill$\blacksquare$

\par

{\bf Proposition 3.4}\quad {\it Let the initial datum
$U^\varepsilon_0$ be a $ \mathcal{F}_0$-measurable random variable
in $L^2(\Omega, \mathcal{H}_\varepsilon)$. Then the solution
$U^\varepsilon(t)$ of the Cauchy problem (\ref{Eq4}) globally exists
in $\mathcal{H}_\varepsilon$, i.e. $\tau^*=+\infty$ almost surely. }
\par
{\bf Proof.}\quad For any given positive $T_0$, consider the case
  that $t<\tau^*\leq T_0$. For any stopping time $\tau$ satisfying
$\tau<\tau^*$, it follows from Proposition 3.3 and the Gronwall
inequality that for arbitrary $t\leq \tau \wedge \tau_R$,
\begin{equation}\label{Eq-3.4-1}
\mathbb{E}\|U_r^\varepsilon(t)\|_{\mathcal{H}_\varepsilon}^2\leq
C(T_0, TrQ_1, TrQ_2, \mathbb{E} \mathcal{E}_r^\varepsilon(0),
\mathbb{E}\|U_r^\varepsilon(0)\|_{\mathcal{H}_\varepsilon}),
\end{equation}
where $\tau_R$ is defined as (\ref{Eq7}).
\par
Moreover, we note that from Frigeri \cite{Frigeri}, for
$r\in(0,\frac{1}{2})$,
$\mathbb{E}\|U_r^\varepsilon\|_{\mathcal{H}_\varepsilon}^2\geq
\frac{1}{2}\mathbb{E}\|U^\varepsilon\|_{\mathcal{H}_\varepsilon}^2$.
Then take $r \in(0,\frac{1}{2})$ sufficiently small such that
(\ref{Eq-3.3-6}) holds. Then for arbitrary $t\leq \tau \wedge
\tau_R$,
\begin{equation}\label{Eq-3.4-2}
\begin{array}{ll}
\mathbb{E}\|U_r^\varepsilon(t)\|_{\mathcal{H}_\varepsilon}^2& \geq
\frac{1}{2}
\mathbb{E}\|U^\varepsilon(t)\|_{\mathcal{H}_\varepsilon}^2\\
&\geq C\mathbb{E}[\|U^\varepsilon(t)\|_{\mathcal{H}_\varepsilon}^2\cdot\chi(\{\tau_R\leq T_0\})]\\
&\geq C \mathbb{E}[\frac{R}{2}\cdot\chi(\{\tau_R\leq T_0\})]\\
&=C\cdot \frac{R}{2}\mathbb{P}\{\tau_R\leq T_0\},
\end{array}
\end{equation}
where $\chi$ is the indicator function.
\par
Therefore, from (\ref{Eq-3.4-1}) and (\ref{Eq-3.4-2}), we  see that
\begin{equation}\label{Eq-3.4-3}
\mathbb{P}\{\tau_R\leq T_0\}\leq \frac{2C(T_0, TrQ_1, TrQ_2,
\mathbb{E} \mathcal{E}_r^\varepsilon(0),
\mathbb{E}\|U^\varepsilon(0)\|_{\mathcal{H}_\varepsilon} )}{CR},
\end{equation}
which implies from the Borel-Cantelli lemma that
\begin{equation}\label{Eq-3.4-4}
\mathbb{P}\{\tau^*\leq T_0\}=0,
\end{equation}
where $\tau^*=\lim\limits_{R\to +\infty}\tau_R$. In other words, we
conclude that
\begin{equation}\label{Eq-3.4-5}
\mathbb{P}\{\tau^*=\infty\}=1.
\end{equation}
Therefore the solution $U^\varepsilon(t)$ of the Cauchy problem
(\ref{Eq4}) globally exists almost surely. This completes the proof.
\hfill$\blacksquare$

\par

{\bf Proposition 3.5} \quad{\it Let the initial datum
$U^\varepsilon_0$ be a $ \mathcal{F}_0$-measurable random variable
in $L^2(\Omega, \mathcal{H}_\varepsilon)$. Then the global solution
$U^\varepsilon(t)$ of the Cauchy problem (\ref{Eq4}) is bounded in
$\mathcal{H}_\varepsilon$ almost surely. }
\par
{\bf Proof.}\quad From Proposition 3.4, we know that the solution
$U^\varepsilon(t)$ of the Cauchy problem (\ref{Eq4}) globally exists
on $[0, +\infty)$ almost surely. Therefore, it follows from
Proposition 3.3 that for arbitrary $t\in[0, +\infty)$,
$$
\frac{d}{dt}\mathbb{E}\|U_r^\varepsilon(t)\|_{\mathcal{H}_\varepsilon}^2+C
\mathbb{E}\|U_r^\varepsilon(t)\|_{\mathcal{H}_\varepsilon}^2\leq
C[TrQ_1 +TrQ_2+1],
$$
which immediately implies from the Gronwall inequality that
\begin{equation}\label{Eq-3.5-1}
\mathbb{E}\|U_r^\varepsilon(t)\|_{\mathcal{H}_\varepsilon}^2\leq
\mathbb{E}\|U_r^\varepsilon(0)\|_{\mathcal{H}_\varepsilon}^2e^{-Ct}+[TrQ_1
+TrQ_2+1](1-e^{-Ct}).
\end{equation}
Note that for $r\in(0,\frac{1}{2})$,
$\mathbb{E}\|U_r^\varepsilon\|_{\mathcal{H}_\varepsilon}^2\geq
\frac{1}{2}\mathbb{E}\|U^\varepsilon\|_{\mathcal{H}_\varepsilon}^2$.
Thus we take $r \in(0,\frac{1}{2})$ sufficiently small such that
(\ref{Eq-3.3-6}) holds. It then follows from (\ref{Eq-3.5-1}) that
Proposition 3.5 holds. \hfill$\blacksquare$
\par

Introduce a space
$$ \Sigma_\varepsilon: =\{U^\varepsilon\in
H_\varepsilon^2(D^\varepsilon)\times
H_\varepsilon^1(D^\varepsilon)\times H^1(\partial
S^\varepsilon)\times H^1(\partial S^\varepsilon)|\; \frac{\partial
u^\varepsilon}{\partial {\bf n}}=\theta^\varepsilon\; \hbox{on} \;
\partial S^\varepsilon \},$$
where $H^2_\varepsilon(D^\varepsilon)$ and
$H^1_\varepsilon(D^\varepsilon)$ denote the space
$H^2(D^\varepsilon)$ and $H^1(D^\varepsilon)$ vanishing on $\partial
D$, respectively.
\par
{\bf Proposition 3.6} \quad{\it Let the initial datum
$U^\varepsilon_0$ be a $ \mathcal{F}_0$-measurable random variable
in $L^2(\Omega, \Sigma_\varepsilon)$. Then the global solution
$U^\varepsilon(t)$ of the Cauchy problem (\ref{Eq4}) is also bounded
in $\Sigma_\varepsilon$ almost surely. }
\par
The proof of Proposition 3.6 is similar as Proposition 3.2,
Proposition 3.3 and Proposition 3.5. It is omitted here.
\par

In the following, for any $T>0$, we consider the solution
$(u^\varepsilon, v^\varepsilon)^T\in L^2(0,T;
H_\varepsilon^1(D^\varepsilon)\times
L_\varepsilon^2(D^\varepsilon))$ of Equation (\ref{Eq2}). Set
$$
\mathcal{X}:=H^1(D)\times L^2(D), \quad \mathcal{Y}:=L^2(D)\times
L^2(D),\quad \mathcal{Z}:=H^{-1}(D)\times L^2(D).
$$
\par
We investigate the behavior of distribution of $(u^\varepsilon,
v^\varepsilon)^T\in L^2(0,T; L_\varepsilon^2(D^\varepsilon)\times
L_\varepsilon^2(D^\varepsilon))$ as $\varepsilon\to 0$, which needs
the tightness of distribution (see \cite{Dudley}). Notice that the
function space changes with $\varepsilon$, which is a difficulty for
obtaining the tightness of distributions. Thus we will treat
$\{\mathcal{L}((u^\varepsilon, v^\varepsilon)^T))\}_{\varepsilon>0}$
as a collection of distributions on $L^2((0,T), \mathcal{Y})$ by
extending $(u^\varepsilon, v^\varepsilon)^T$ to the whole domain
$D$, whose distribution is defined as
$\mathcal{L}((\tilde{u}^\varepsilon,
\tilde{v}^\varepsilon)^T)(A)=\mathbb{P}\{\omega:
(\tilde{u}^\varepsilon(\cdot,\cdot, \omega),
\tilde{v}^\varepsilon(\cdot,\cdot, \omega))^T\in A\}$ for the Borel
set $A\in L^2((0,T), \mathcal{Y})$.
\par
{\bf Proposition 3.7 (Tightness of distribution)}  \quad{\it Let the
initial datum $U^\varepsilon_0$ be a $ \mathcal{F}_0$-measurable
random variable in $L^2(\Omega, \mathcal{H}_\varepsilon)$, which is
independent of $W(t)$ with
$\mathbb{E}\|U^\varepsilon_0\|_{\mathcal{H}_\varepsilon}^2<\infty$.
Then for any $T>0$, $\mathcal{L}((u^\varepsilon, v^\varepsilon)^T)$,
the distribution of $(u^\varepsilon, v^\varepsilon)^T$, is tight in
$L^2((0,T), \mathcal{Y})\bigcap C((0,T),\mathcal{Z})$.}
\par
{\bf Proof.}\quad Firstly, we claim that $(u^\varepsilon,
v^\varepsilon)^T$ is bounded almost surely in
$$
G:=L^2(0,T; \mathcal{X})\bigcap(W^{1,2}(0,T; \mathcal{Z})+W^{\alpha,
4}(0,T; \mathcal{Y})),
$$
where $W^{1,2}(0,T; \mathcal{Z})$ is a Banach space endowed with the norm
$$
\|\varphi\|_{W^{1,2}(0,T; \mathcal{Z})}^2=\|\varphi\|_{L^2(0,T;
\mathcal{Z})}^2+\|\frac{d\varphi}{dt}\|_{L^2(0,T;
\mathcal{Z})}^2<\infty, \quad \forall\quad \varphi\in W^{1,2}(0,T;
\mathcal{Z}),
$$
and $W^{\alpha,4}(0,T; \mathcal{Y})$ is another Banach space with
$\alpha\in (\frac{1}{4}, \frac{1}{2})$ endowed with the norm
$$
\|\varphi\|_{W^{\alpha,4}(0,T; \mathcal{Y})}^4=\|\varphi\|_{L^4(0,T;
\mathcal{Y})}^4+\int_0^T\int_0^T\frac{\|\varphi(t)-\varphi(s)\|_\mathcal{Y}^4}{|t-s|^{1+4\alpha}}dsdt
<\infty, \quad \forall \quad \varphi\in W^{\alpha,4}(0,T;
\mathcal{Y}).
$$
\par
By Proposition 3.5, we know that $(u^\varepsilon,
v^\varepsilon)^T$ is bounded in $L^2(0,T; \mathcal{X})$ almost
surely. Therefore, in the following, we only need to prove that
$(u^\varepsilon, v^\varepsilon)^T$ is bounded in $W^{1,2}(0,T;
\mathcal{Z})+W^{\alpha, 4}(0,T; \mathcal{Y})$ almost surely.
\par

Denote by $P$ the projection operator from $U^\varepsilon$ to
$(u^\varepsilon, v^\varepsilon)^T$, i.e.,
$PU^\varepsilon=(u^\varepsilon, v^\varepsilon)^T$. Write Equation
(\ref{Eq4}) as
$$
U^\varepsilon(t)=U^\varepsilon(0)+\int_0^tA^\varepsilon
U^\varepsilon(\tau)d\tau+\int_0^tF^\varepsilon(U^\varepsilon(\tau))d\tau+\int_0^tdW(\tau).
$$
Then
\begin{equation}\label{Eq-3.7-1}
PU^\varepsilon(t)=PU^\varepsilon(0)+\int_0^t[PA^\varepsilon
U^\varepsilon(\tau)+PF^\varepsilon(U^\varepsilon(\tau))]d\tau+\int_0^tPdW(\tau).
\end{equation}
\par
Denote
\begin{equation}\label{Eq-3.7-2}
I_1:=\int_0^t[PA^\varepsilon
U^\varepsilon(\tau)+PF^\varepsilon(U^\varepsilon(\tau))]d\tau,
\end{equation}
and
\begin{equation}\label{Eq-3.7-3}
I_2:= \int_0^tPdW(\tau).
\end{equation}
\par
For $I_1$, it follows from Proposition 3.5 and Proposition 3.6 that
\begin{equation}\label{Eq-3.7-4}
\begin{array}{ll}
&\mathbb{E}\|I_1\|_{L^{2}(0,T;
\mathcal{Z})}^2\\
=&\mathbb{E}\int_0^T\|I_1(\tau)\|_\mathcal{Z}^2d\tau\\
=&\mathbb{E}\int_0^T\|\int_0^\tau[PA^\varepsilon
U^\varepsilon(s)+PF^\varepsilon(U^\varepsilon(s))]ds\|_\mathcal{Z}^2 d\tau\\
\leq &\mathbb{E}\int_0^T[\|\int_0^\tau
v^\varepsilon(s)ds\|_{H^{-1}(D)}^2+\|\int_0^\tau [\bigtriangleup
u^\varepsilon(s)-u^\varepsilon(s)-v^\varepsilon(s)+\sin
u^\varepsilon(s)]ds\|_{L^2(D)}^2]d\tau\\
\leq& \mathbb{E}\int_0^T[\int_0^\tau
\|v^\varepsilon(s)\|_{L^2(D)}^2ds+\int_0^\tau[\|\bigtriangleup
u^\varepsilon(s)\|_{L^2(D)}^2+\|u^\varepsilon(s)\|_{L^2(D)}^2\\
&+\|v^\varepsilon(s)\|_{L^2(D)}^2
+\|u^\varepsilon(s)\|_{L^2(D)}^2]ds]d\tau\\
\leq& C_T,
\end{array}
\end{equation}
and
\begin{equation}\label{Eq-3.7-5}
\begin{array}{ll}
&\mathbb{E}\|\frac{d I_1}{dt}\|_{L^{2}(0,T;
\mathcal{Z})}^2\\
=&\mathbb{E}\int_0^T\|\frac{d I_1}{d\tau}\|_\mathcal{Z}^2d\tau\\
=&\mathbb{E}\int_0^T\|PA^\varepsilon
U^\varepsilon(\tau)+PF^\varepsilon(U^\varepsilon(\tau))\|_\mathcal{Z}^2 d\tau\\
\leq& \mathbb{E}\int_0^T[\|
v^\varepsilon(\tau)\|_{H^{-1}(D)}^2+\|\bigtriangleup
u^\varepsilon(\tau)-u^\varepsilon(\tau)-v^\varepsilon(\tau)+\sin
u^\varepsilon(\tau)\|_{L^2(D)}^2]d\tau\\
\leq &\mathbb{E}\int_0^T[
\|v^\varepsilon(\tau)\|_{L^2(D)}^2+\|\bigtriangleup
u^\varepsilon(\tau)\|_{L^2(D)}^2+\|u^\varepsilon(\tau)\|_{L^2(D)}^2\\
&+\|v^\varepsilon(\tau)\|_{L^2(D)}^2
+\|u^\varepsilon(\tau)\|_{L^2(D)}^2]d\tau\\
\leq &C_T.
\end{array}
\end{equation}
Here and hereafter, $C_T$ denotes various positive constants
depending on the given $T>0$. Then combining (\ref{Eq-3.7-4}) and
(\ref{Eq-3.7-5}), we deduce that
\begin{equation}\label{Eq-3.7-6}
\mathbb{E}\|I_1\|_{W^{1,2}(0,T;
\mathcal{Z})}^2=\mathbb{E}\|I_1\|_{L^{2}(0,T;
\mathcal{Z})}^2+\mathbb{E}\|\frac{d I_1}{dt}\|_{L^{2}(0,T;
\mathcal{Z})}^2\leq C_T.
\end{equation}
\par
Now we consider $I_2$. Put $M^\varepsilon(s,
t)=\int_s^tdW(\tau)$. Then using the Burkholder-Davis-Gundy
inequality and the H$\ddot{o}$lder inequality, we have
\begin{equation}\label{Eq-3.7-7}
\begin{array}{ll}
\mathbb{E}\|PM^\varepsilon(s, t)\|_{\mathcal{Y}}^4&= \mathbb{E}\|\int_s^tPdW(\tau)\|_{\mathcal{Y}}^4\\
&=\mathbb{E}\|\int_s^tdW_1(\tau)\|_{L^2(D)}^4\\
& \leq C \mathbb{E}(\int_s^tTrQ_1d\tau)^2\\
& \leq C \mathbb{E}(\int_s^t1^2d\tau\cdot\int_s^t(TrQ_1)^2d\tau)\\
& \leq C(t-s)^2.
\end{array}
\end{equation}
Thus, it follows from (\ref{Eq-3.7-7}) that
\begin{equation}\label{Eq-3.7-8}
\begin{array}{ll}
\mathbb{E}\|I_2\|_{L^4(0,T;
\mathcal{Y})}^4&=\mathbb{E}\int_0^T\|I_2\|_{\mathcal{Y}}^4dt\\
&=\mathbb{E}\int_0^T\|\int_0^tPdW(\tau)\|_{\mathcal{Y}}^4dt\\
&=\mathbb{E}\int_0^T\|PM^\varepsilon(0,t)\|_{\mathcal{Y}}^4dt\\
&\leq \int_0^TC(t-0)^2dt\\
&\leq C_T.
\end{array}
\end{equation}
Also, for $\alpha\in (\frac{1}{4}, \frac{1}{2})$, by
(\ref{Eq-3.7-7}), we have
\begin{equation}\label{Eq-3.7-9}
\begin{array}{ll}
\mathbb{E}\int_0^T\int_0^T\frac{\|I_2(t)-I_2(s)\|_\mathcal{Y}^4}{|t-s|^{1+4\alpha}}dsdt
&=
\mathbb{E}\int_0^T\int_0^T\frac{\|PM^\varepsilon(0,t)-PM^\varepsilon(0,s)\|_\mathcal{Y}^4}{|t-s|^{1+4\alpha}}dsdt\\
&=
\mathbb{E}\int_0^T\int_0^T\frac{\|PM^\varepsilon(s,t)\|_\mathcal{Y}^4}{|t-s|^{1+4\alpha}}dsdt\\
&\leq \int_0^T\int_0^T\frac{C(t-s)^2}{(t-s)^{1+4\alpha}}dsdt\\
&\leq C\int_0^T\int_0^T (t-s)^{1-4\alpha}dsdt\\
&\leq C_T.
\end{array}
\end{equation}
Therefore, it follows from (\ref{Eq-3.7-8}) and (\ref{Eq-3.7-9}) that for
arbitrary $\alpha\in (\frac{1}{4},\frac{1}{2})$,
\begin{equation}\label{Eq-3.7-10}
\mathbb{E}\|I_2\|_{W^{\alpha,4}(0,T;
\mathcal{Y})}^4=\mathbb{E}\|I_2\|_{L^4(0,T;
\mathcal{Y})}^4+\mathbb{E}\int_0^T\int_0^T\frac{\|I_2(t)-I_2(s)\|_\mathcal{Y}^4}{|t-s|^{1+4\alpha}}dsdt<C_T.
\end{equation}
Immediately from (\ref{Eq-3.7-1})-(\ref{Eq-3.7-3}),
(\ref{Eq-3.7-6})and (\ref{Eq-3.7-10}), we obtain that
$(u^\varepsilon, v^\varepsilon)^T$ is bounded in $W^{1,2}(0,T;
\mathcal{Z})+W^{\alpha, 4}(0,T; \mathcal{Y})$ almost surely, which
completes the verification of the claim that $(u^\varepsilon,
v^\varepsilon)^T$ is bounded almost surely in $G=L^2(0,T;
\mathcal{X})\bigcap(W^{1,2}(0,T; \mathcal{Z})+W^{\alpha, 4}(0,T;
\mathcal{Y}))$.
\par
By the Chebyshev inequality, we see that for any $\rho>0$,
there exists a bounded set $K_\rho\subset G$ such that
$\mathbb{P}\{(u^\varepsilon, v^\varepsilon)^T\in K_\rho\}>1-\rho$.
Moreover, notice that
$$
L^2(0,T; \mathcal{X)}\bigcap W^{1,2}(0,T; \mathcal{Z})\subset L^2(0,
T; \mathcal{Y})\bigcap C(0,T; \mathcal{Z}),
$$
and for $\alpha\in (\frac{1}{4}, \frac{1}{2})$,
$$
L^2(0,T; \mathcal{X})\bigcap W^{\alpha,4}(0,T; \mathcal{Y})\subset
L^2(0, T; \mathcal{Y})\bigcap C(0,T; \mathcal{Z}).
$$
We conclude that $K_\rho$ is compact in $L^2(0, T;
\mathcal{Y})\bigcap C(0,T; \mathcal{Z})$. Thus
$\mathcal{L}((u^\varepsilon, v^\varepsilon)^T)$ is tight in $L^2(0,
T; \mathcal{Y})\bigcap C(0,T; \mathcal{Z})$. \hfill$\blacksquare$
\par

\renewcommand{\theequation}{\thesection.\arabic{equation}}
\setcounter{equation}{0}

\section{Effective Model}

\quad\quad In this section, we will use the two-scale method to
derive the effective  homogenized   equation of Equation
(\ref{Eq1}), in the sense of probability distribution. The solutions
of the microscopic model Equation (\ref{Eq1}) converge to those of
the effective homogenized   equation in probability
distribution, as the size of small holes $\varepsilon$ diminishes to
zero. The main result is as follows.
\par
{\bf Theorem 4.1 (Homogenized   model)}\quad {\it Let
$(u^\varepsilon, \delta^\varepsilon)^T$ be the solution of Equation
(\ref{Eq1}). Then for any $T>0$, the distribution
$\mathcal{L}(\tilde{u}^\varepsilon)$ converges weakly to $\mu$ in
$L^2(0, T; L^2(D))$ as $\varepsilon\to 0$, with $\mu$ being the
distribution of the solution $V$  of the following
  homogenizied   equation
\begin{equation}\label{Eq-E}
\left\{
\begin{array}{l}
V_{tt}(t,x)+ V_t(t,x)-\nu^{-1}div_x\mathcal{A}^*(\bigtriangledown_x
V(t,x)) + V(t,x)-\sin(V(t,x)) =\nu \dot{W}_1,\quad \hbox{on}\quad D\\
V(t,x)=0,\quad \hbox{on}\quad \partial D\\
V(0,x)=\frac{u_0}{\nu},\quad V_t(0,x)=\frac{v_0}{\nu},\quad
\hbox{on}\quad D,
\end{array}
\right.
\end{equation}
where the effective matrix $\mathcal{A}^*=(A_{ij}^*)$ given by
(\ref{Eq-A}), $u_0$ and $v_0$ are the initial data supplemented in
Equation (\ref{Eq3}), and the constant $\nu=\frac{|Y^*|}{|Y|}$ with
$|Y|$ and $|Y^*|$ the Lebesgue measure of $Y$ and $Y^*$
respectively. }
\par
In the following, we will prove Theorem 4.1. We first
provide some preliminaries. We will denote by $C_{per}^\infty(Y)$ the
space of infinitely differentiable functions in $\mathbb{R}^3$ that
are periodic in $Y$. We also denote $L_{per}^2(Y)$ or $H_{per}^1(Y)$
the completion of $C_{per}^\infty(Y)$ in the usual norm of $L^2(Y)$
or $H^1(Y)$, respectively. In addition, we denote $D_T=[0,T]\times
D$.
\par
{\bf Definition 4.1}$^{\cite{A}}$\quad {\it A sequence of functions
$u^\varepsilon(t,x)$ in $L^2(D_T)$ is called to be two-scale
convergent to a limit $u(t,x, y)\in L^2(D_T\times Y)$, if for any
function $\varphi(x, y)\in C_0^\infty(D_T, C_{per}^\infty)$,
$$ \lim\limits_{\varepsilon\to
0}\int_{D_T}u^\varepsilon(t,x)\varphi(t,x,\frac{x}{\varepsilon})dxdt=\frac{1}{|Y|}\int_{D_T}\int_Y
u(t,x,y)\varphi(t,x,y)dydxdt,
$$
which is denoted by $u^\varepsilon \xlongrightarrow{2-s} u$.}
\par
{\bf Lemma 4.1}$^{\cite{A}}$\quad{\it Let $u^\varepsilon$ be a
bounded sequence in $L^2(D_T)$. Then there exists a function $u\in
L^2(D_T\times Y)$ and a subsequence $u_{\varepsilon_k}\to 0$ as
$k\to \infty$ such that $u_{\varepsilon_k}\xlongrightarrow{2-s} u$.}
\par
{\bf Lemma 4.2}$^{\cite{A}}$\quad{\it If $u^\varepsilon
\xlongrightarrow{2-s} u$, then $u^\varepsilon \rightharpoonup
\overline{u}(t,x)=\frac{1}{|Y|}\int_Yu(t,x,y)dy$.}
\par
{\bf Lemma 4.3}$^{\cite{A}}$\quad {\it Let $v^\varepsilon$ be a
sequence in $L^2(D_T)$ that two-scale converges to a limit
$v(x,y)\in L^2(D_T\times Y)$. Further assume that
$$
\lim\limits_{\varepsilon\to
0}\int_{D_T}|v^\varepsilon(t,x)|^2dxdt=\frac{1}{|Y|}\int_{D_T}\int_Y|v(t,x,y)|^2dydxdt.
$$
Then for any sequence $u^\varepsilon\in L^2(D_T)$, which two-scale
converges to a limit $u\in L^2(D_T\times Y)$, we have
$$
u^\varepsilon v^\varepsilon \rightharpoonup
\frac{1}{|Y|}\int_Yu(\cdot,\cdot, y)v(\cdot,\cdot, y)dy,\;
\hbox{as}\; \varepsilon\to 0\; \hbox{in} \;L^2(D_T).$$ }
\par
{\bf Lemma 4.4}$^{\cite{A}}$\quad {\it Let $u^\varepsilon$ be a
sequence of functions defined on $[0, T]\times D^\varepsilon$ which
is bounded in $L^2(0, T; H_{\varepsilon}^1(D^\varepsilon))$. There
exists $u(t,x)\in H_0^1(D_T)$, $u_1(t,x,y)\in L^2(D_T;H_{per}^1(Y))$
and a subsequence $u^{\varepsilon_k}$ with $\varepsilon_k\to 0$ as
$k\to \infty$, such that
$$
\tilde{u}^{\varepsilon_k}(t,x)  \xlongrightarrow{2-s} \chi(Y)u(t,x),
\quad k\to \infty,
$$
and
$$
\widetilde{\bigtriangledown_xu^{\varepsilon_k}}
\xlongrightarrow{2-s} \chi(Y)[\bigtriangledown_x
u(t,x)+\bigtriangledown_y u_1(t,x,y)],\quad k\to \infty,
$$
where $\chi(Y)$ is the indicator function as defined in Section 2. }
\par
For $h\in H^{-1/2}(\partial S)$ and $Y$-periodic, define
$\lambda_h:=\frac{1}{|Y|}\int_{\partial S} h(x) dx$. Also, for $h\in
L^2(\partial S)$ and $Y$-periodic, define $\lambda_h^\varepsilon \in
H^{-1}(D)$ as $\langle\lambda_h^\varepsilon, \varphi\rangle_{H^{-1},
H_0^1}=\varepsilon\int_{\partial
S^\varepsilon}h(\frac{x}{\varepsilon})\varphi(x)dx$ with any
$\varphi\in H_0^1(D)$.
\par
{\bf Lemma 4.5}$^{\cite{WD}}$\quad {\it Let $\varphi^\varepsilon$ be
a sequence in $H_0^1(D)$ such that $\varphi^\varepsilon
\rightharpoonup \varphi$ in $H_0^1(D)$ as $\varepsilon\to 0$. Then
$$
\langle\lambda_h^\varepsilon,
\varphi^\varepsilon|_{D^\varepsilon}\rangle\longrightarrow
\lambda_h\int_D\varphi dx,\quad \hbox{as}\quad \varepsilon \to 0.
$$
 }
\par
{\bf Lemma 4.6 (Prohorov Theorem)}$^{\cite{DZ1}}$\quad {\it Suppose
$\mathcal{M}$ is a separable Banach space. The set of probability
measures $\{\mathcal{L}(X_n)\}_n$ on $(\mathcal{M},
\mathcal{B}(\mathcal{M}))$ is relatively compact if and only if
$\{X_n\}$ is tight.}
\par
{\bf Lemma 4.7 (Skorohod Theorem)}$^{\cite{DZ1}}$\quad {\it For an
arbitrary sequence of Probability measures $\{\mu_n\}$ on
$(\mathcal{M}, \mathcal{B}(\mathcal{M}))$ weakly converges to
probability measures $\mu$, there exists a probability space
$(\Omega, \mathcal{F}, \mathbb{P})$ and random variables, $X$,
$X_1$, $X_2$, $\cdots$, $X_n$, $\cdots$ such that $X_n$ distributes
as $\mu_n$ and $X$ distributes as $\mu$ and $\lim\limits_{n\to
\infty}X_n=X$, $\mathbb{P}$-a.s.}
\par

\medskip

\vspace{0.7cm}

{\bf Proof of   Theorem 4.1.}
\par

\par
Let $(u^\varepsilon, \delta^\varepsilon)^T$ be the
solution of Equation (\ref{Eq1}). On the one hand, as in \cite{WD}, by the proof of Proposition
3.7, for any $\rho>0$, there is a bounded set $K_\rho\subset G$
which is compact in $L^2(0, T; \mathcal{Y})$ such that
$\mathbb{P}\{(\tilde{u}^\varepsilon,\tilde{v}^\varepsilon)^T\in
K_\rho\}>1-\rho$. According to Lemma 4.6 and Lemma 4.7, we know that
for any sequence $\{\varepsilon_j\}_{j=1}^{j=\infty}$ with
$\varepsilon_j\to 0$ as $j\to \infty$, there exists a subsequence
$\{\varepsilon_{j(k)}\}_{k=1}^{k=\infty}$, random variables
$\{(\tilde{u}^{\varepsilon_{j(k)}}_*,
\tilde{v}^{\varepsilon_{j(k)}}_*)^T\}\subset L^2(0, T;
L_\varepsilon^2(D^\varepsilon)\times
L_\varepsilon^2(D^\varepsilon))$ and $(u_*,v_*)^T\in L^2(0,T;
\mathcal{Y})$ defined on a new probability space $(\Omega_*,
\mathcal{F}_*, \mathbb{P}_*)$, such that for almost all $\omega\in
\Omega_*$,
$$
\mathcal{L}((\tilde{u}^{\varepsilon_{j(k)}}_*,
\tilde{v}^{\varepsilon_{j(k)}}_*)^T)=
\mathcal{L}((\tilde{u}^{\varepsilon_{j(k)}},
\tilde{v}^{\varepsilon_{j(k)}})^T),
$$
and
\begin{equation}\label{sc}
(\tilde{u}^{\varepsilon_{j(k)}}_*,
\tilde{v}^{\varepsilon_{j(k)}}_*)^T\longrightarrow (u_*,v_*)^T\quad
\hbox{in}\quad L^2(0,T; \mathcal{Y}) \quad \hbox{as}\quad k\to
\infty.
\end{equation}
In the meantime, $(\tilde{u}^{\varepsilon_{j(k)}}_*,
\tilde{v}^{\varepsilon_{j(k)}}_*)^T $   solves
$$
\left\{
\begin{array}{l}
dPU^\varepsilon=PA^\varepsilon U^\varepsilon dt+PF^\varepsilon(U^\varepsilon)dt+ PdW(t),\\
PU^\varepsilon(0)=PU_0,
\end{array}
\right.
$$
with $W$ being replaced by a Wiener process $W_*$, defined on the
probability space $(\Omega_*, \mathcal{F}_*, \mathbb{P}_*)$ but with the
same distributions as $W$. Here $P$ is the projection operator from
$U^\varepsilon$ to $(u^\varepsilon, v^\varepsilon)^T$ as defined in
the proof of Proposition 3.7.
\par
On the other hand, for $u^\varepsilon$ in the set $K_\rho$, it
follows from Lemma 4.1 and Lemma 4.4 that  there exist $u(t,x)\in
H_0^1(D_T)$ and $u_1(t,x,y)\in L^2(D_T, H^1_{per}(Y))$ such that
\begin{equation}\label{2s1}
\tilde{u}^{\varepsilon_k}(t,x)  \xlongrightarrow{2-s} \chi(Y)u(t,x),
\quad k\to \infty,
\end{equation}
and
\begin{equation}\label{2s2}
\widetilde{\bigtriangledown_xu^{\varepsilon_k}}
\xlongrightarrow{2-s} \chi(Y)[\bigtriangledown_x
u(t,x)+\bigtriangledown_y u_1(t,x,y)],\quad k\to \infty.
\end{equation}
Furthermore, from Lemma 4.2, it follows that
$$
\tilde{u}^{\varepsilon_k}(t,x) \rightharpoonup
\frac{1}{|Y|}\int_Y\chi(Y)u(t,x)dy=\frac{1}{|Y|}\int_Y\chi(Y)dy\cdot
u(t,x)= \nu u(t,x),\quad\hbox{in}\quad L^2(D_T),
$$
which from the compactness of $K_\rho$ immediately implies that
\begin{equation}\label{s-c}
\tilde{u}^{\varepsilon_k}(t,x) \longrightarrow  \nu
u(t,x),\quad\hbox{in}\quad L^2(D_T).
\end{equation}
Then combining the relationship of $u^\varepsilon$ and
$v^\varepsilon$, (\ref{sc}) and (\ref{s-c}), we have
\begin{equation}\label{uv}
u_*=\nu u \quad \hbox{and} \quad v_*=\nu u_t.
\end{equation}
\par
Now, in the probability space $(\Omega, \mathcal{F}, \mathbb{P})$,
we put $\Omega_\rho=\{\omega\in\Omega:
\tilde{u}^\varepsilon(\omega)\in K_\rho\}$,
$\mathcal{F}_\rho=\{F\bigcap\Omega_\rho: F\in \mathcal{F}\}$, and
$\mathbb{P}_\rho(F)=\frac{\mathbb{P}(F\bigcap\Omega_\rho)}{\mathbb{P}(\Omega_\rho)}$,
for $F\in \mathcal{F}_\rho$. Then $(\Omega_\rho, \mathcal{F}_\rho,
\mathbb{P}_\rho)$ forms  a new probability space, whose expectation
operator is denoted by $\mathbb{E}_\rho$. In the following, we will
work in the probability space $(\Omega_\rho, \mathcal{F}_\rho,
\mathbb{P}_\rho)$ in stead of $(\Omega, \mathcal{F}, \mathbb{P})$.
\par
In Equation (\ref{E-v}), we choose the test function $\varphi$ as
$\varphi^\varepsilon(t,x)=\phi(t,x)+\varepsilon\Phi(t,x,\frac{x}{\varepsilon})$
with $\phi(t,x)\in C_0^\infty(D_T)$ and $\Phi(t,x,y)\in
C_0^\infty(D_T; C_{per}^\infty(Y))$. Also, we notice that
(\ref{s-c}) and $\chi(D^\varepsilon)\rightharpoonup \nu$ in
$L^\infty(D)$. Then we have
\begin{equation}\label{Eq-4.1-1}
\begin{array}{lll}
\int_0^T\int_{D^\varepsilon}u^\varepsilon_{tt}\varphi^\varepsilon
dxdt&=&\int_0^T\int_{D^\varepsilon}u^\varepsilon_{tt}[\phi(t,x)+\varepsilon\Phi(t,x,\frac{x}{\varepsilon})]
dxdt\\
&=&-\int_0^T\int_{D^\varepsilon}[u^\varepsilon_{t}\phi_t(t,x)+\varepsilon
u^\varepsilon_{t}\Phi_t(t,x,\frac{x}{\varepsilon}) ]dxdt\\
&=&\int_0^T\int_{D^\varepsilon}[u^\varepsilon\phi_{tt}(t,x)+\varepsilon
u^\varepsilon\Phi_{tt}(t,x,\frac{x}{\varepsilon})] dxdt\\
&=&\int_0^T\int_{D}[\tilde{u}^\varepsilon\phi_{tt}(t,x)+\varepsilon
\tilde{u}^\varepsilon\Phi_{tt}(t,x,\frac{x}{\varepsilon}) ]dxdt\\
&\to& \int_0^T\int_{D}\nu u(t,x) \phi_{tt}(t,x)dxdt\\
&=&\int_0^T\int_{D}\nu u_{tt}(t,x) \phi(t,x)dxdt, \quad
\hbox{as}\quad \varepsilon\to 0,
\end{array}
\end{equation}

\begin{equation}\label{Eq-4.1-2}
\begin{array}{lll}
\int_0^{T}\int_{D^\varepsilon}u^\varepsilon_{t}\varphi^\varepsilon
dxdt &=&
\int_0^{T}\int_{D^\varepsilon}u^\varepsilon_{t}[\phi(t,x)+\varepsilon\Phi(t,x,\frac{x}{\varepsilon})]
dxdt\\
&=&-\int_0^{T}\int_{D^\varepsilon}[u^\varepsilon\phi_t(t,x)+\varepsilon
u^\varepsilon\Phi_t(t,x,\frac{x}{\varepsilon}) ]dxdt\\
&=&-\int_0^{T}\int_{D}[\tilde{u}^\varepsilon\phi_t(t,x)+\varepsilon
\tilde{u}^\varepsilon\Phi_t(t,x,\frac{x}{\varepsilon}) ]dxdt\\
&\to&
-\int_0^{T}\int_{D}\nu u(t,x) \phi_t(t,x)dxdt\\
&=& \int_0^{T}\int_{D}\nu u_t(t,x) \phi(t,x)dxdt, \quad
\hbox{as}\quad \varepsilon\to 0,
\end{array}
\end{equation}
and
\begin{equation}\label{Eq-4.1-3}
\begin{array}{lll}
\int_0^{T}\int_{D^\varepsilon}u^\varepsilon\varphi^\varepsilon dxdt
&=&
\int_0^{T}\int_{D^\varepsilon}u^\varepsilon[\phi(t,x)+\varepsilon\Phi(t,x,\frac{x}{\varepsilon})]
dxdt\\
&=&
\int_0^{T}\int_{D}[\tilde{u}^\varepsilon\phi(t,x)+\varepsilon\tilde{u}^\varepsilon\Phi(t,x,\frac{x}{\varepsilon})]
dxdt\\
&\to& \int_0^{T}\int_{D}\nu u(t,x)\phi(t,x)dxdt, \quad
\hbox{as}\quad \varepsilon\to 0.
\end{array}
\end{equation}
\par
For $\varphi^\varepsilon$, we have
\begin{equation}\label{Eq-4.1-4}
\begin{array}{lcl}
\bigtriangledown_x\varphi^\varepsilon&=&\bigtriangledown_x\phi(t,x)
+\varepsilon\bigtriangledown_x\Phi(t,x,\frac{x}{\varepsilon})\\
&=& \bigtriangledown_x\phi(t,x)+
\bigtriangledown_y\Phi(t,x,y)\\
& \xlongrightarrow{2-s}& \bigtriangledown_x\phi(t,x)+
\bigtriangledown_y\Phi(t,x,y), \quad \hbox{as}\quad \varepsilon\to
0,
\end{array}
\end{equation}
and
\begin{equation}\label{Eq-4.1-5}
\begin{array}{ll}
\lim\limits_{\varepsilon\to
0}\int_{D_T}|\bigtriangledown_x\varphi^\varepsilon|^2dxdt&=\lim\limits_{\varepsilon\to
0}\int_{D_T}|\bigtriangledown_x\phi(t,x)
+\varepsilon\bigtriangledown_x\Phi(t,x,\frac{x}{\varepsilon})|^2dxdt\\
&=\frac{1}{|Y|}\int_{D_T}\int_Y|\bigtriangledown_x\phi(t,x)+
\bigtriangledown_y\Phi(t,x,y)|^2dydxdt.
\end{array}
\end{equation}
Then it follows from (\ref{2s2}), (\ref{Eq-4.1-4}),(\ref{Eq-4.1-5})
and Lemma 4.3 that
\begin{equation}\label{Eq-4.1-6}
\begin{array}{ll}
&\int_0^{T}\int_{D^\varepsilon}\bigtriangledown u^\varepsilon
\bigtriangledown\varphi^\varepsilon dxdt\\
=& \int_0^{T}\int_{D^\varepsilon}\bigtriangledown_x u^\varepsilon
\bigtriangledown_x[\phi(t,x)+\varepsilon\Phi(t,x,\frac{x}{\varepsilon})]
dxdt\\
=& \int_0^{T}\int_{D}\widetilde{\bigtriangledown_x u^\varepsilon}
[\bigtriangledown_x\phi(t,x)+ \bigtriangledown_y\Phi(t,x,y)]
dxdt\\
\to& \int_0^{T}\int_{D}\frac{1}{|Y|}\int_Y
\chi(Y)[\bigtriangledown_x u(t,x)+\bigtriangledown_y
u_1(t,x,y)][\bigtriangledown_x\phi(t,x)+
\bigtriangledown_y\Phi(t,x,y)]dydxdt\\
=& \frac{1}{|Y|} \int_0^{T}\int_{D}\int_{Y^*} [\bigtriangledown_x
u(t,x)+\bigtriangledown_y u_1(t,x,y)][\bigtriangledown_x\phi(t,x)+
\bigtriangledown_y\Phi(t,x,y)]dydxdt, \quad \hbox{as}\quad
\varepsilon\to 0.
\end{array}
\end{equation}
\par
From (\ref{s-c}) and note that $\sin u$ is continuous and satisfies the
global Lipshitz condition with respective to $u$, we have
\begin{equation}\label{Eq-4.1-7}
\begin{array}{lll}
\int_0^{T}\int_{D^\varepsilon}\sin u^\varepsilon\varphi^\varepsilon
dxdt &=& \int_0^{T}\int_{D^\varepsilon}\sin
u^\varepsilon[\phi(t,x)+\varepsilon\Phi(t,x,\frac{x}{\varepsilon})]
dxdt\\
&=&
\int_0^{T}\int_{D}[\sin\tilde{u}^\varepsilon\phi(t,x)+\varepsilon
\sin\tilde{u}^\varepsilon\Phi(t,x,\frac{x}{\varepsilon})]
dxdt\\
&\to& \int_0^{T}\int_{D}\sin(\nu u(t,x))\phi(t,x)dxdt, \quad
\hbox{as}\quad \varepsilon\to 0.
\end{array}
\end{equation}
\par
Also realize that
\begin{equation}\label{Eq-4.1-8}
\begin{array}{lll}
\int_0^{T}\int_{D^\varepsilon}\dot{W}_1\varphi^\varepsilon dxdt &=&
\int_0^{T}\int_{D^\varepsilon}\dot{W}_1[\phi(t,x)+\varepsilon\Phi(t,x,\frac{x}{\varepsilon})]
dxdt\\
&=&
\int_0^{T}\int_{D}[\chi(D^\varepsilon)\phi(t,x)+\varepsilon\chi(D^\varepsilon)\Phi(t,x,\frac{x}{\varepsilon})]
dxdW_1(t)\\
&\to& \int_0^{T}\int_{D}\nu \varphi(t,x)dxdW_1(t), \quad
\hbox{as}\quad \varepsilon\to 0.
\end{array}
\end{equation}
\par
Moreover, from Proposition 3.6 and Lemma 4.5, we have
\begin{equation}\label{Eq-4.1-9}
\begin{array}{lll}
\varepsilon^2\int_0^{T}\int_{\partial
S^\varepsilon}\delta^\varepsilon_{tt}\varphi^\varepsilon dxdt &=&
-\varepsilon^2\int_0^{T}\int_{\partial
S^\varepsilon}\delta^\varepsilon_{t}\varphi_t^\varepsilon dxdt\\
&=& -\varepsilon\int_{\partial
S^\varepsilon}\varepsilon\cdot\int_0^{T}\theta^\varepsilon\varphi_t^\varepsilon
dtdx\\
&=& -\varepsilon\langle\lambda_1^\varepsilon,
\int_0^{T}\theta^\varepsilon\varphi_t^\varepsilon
dt|_{D^\varepsilon} \rangle\\
&\to& 0, \quad \hbox{as}\quad \varepsilon\to 0,
\end{array}
\end{equation}

\begin{equation}\label{Eq-4.1-10}
\begin{array}{lll}
\varepsilon^2\int_0^{T}\int_{\partial
S^\varepsilon}\delta^\varepsilon\varphi^\varepsilon dxdt &=&
\varepsilon\int_{\partial
S^\varepsilon}\varepsilon\cdot\int_0^{T}\delta^\varepsilon\varphi^\varepsilon
dtdx\\
&=& \varepsilon\langle\lambda_1^\varepsilon,
\int_0^{T}\delta^\varepsilon\varphi^\varepsilon dt|_{D^\varepsilon}\rangle\\
&\to& 0, \quad \hbox{as}\quad \varepsilon\to 0,
\end{array}
\end{equation}

\begin{equation}\label{Eq-4.1-11}
\begin{array}{lll}
\varepsilon^2\int_0^{T}\int_{\partial
S^\varepsilon}u^\varepsilon_t\varphi^\varepsilon dxdt&=&
-\varepsilon\int_{\partial S^\varepsilon}\varepsilon\cdot
\int_0^{T}u^\varepsilon\varphi_t^\varepsilon dt dx\\
&=& -\varepsilon\langle \lambda_1^\varepsilon,
\int_0^{T}u^\varepsilon \varphi_t^\varepsilon
dt|_{D^\varepsilon}\rangle\\
&\to& 0, \quad \hbox{as}\quad \varepsilon\to 0,
\end{array}
\end{equation}
and
\begin{equation}\label{Eq-4.1-12}
\begin{array}{lll}
\varepsilon^2\int_0^{T}\int_{\partial
S^\varepsilon}\dot{W}_2\varphi^\varepsilon dxdt &=&
\varepsilon\int_{\partial S^\varepsilon} \varepsilon\cdot\int_0^{T}
\varphi^\varepsilon dW_2(t) dx\\
&=& \varepsilon\langle \lambda_1^\varepsilon,\int_0^{T}
\varphi^\varepsilon dW_2(t)|_{D^\varepsilon}\rangle\\
&\to& 0, \quad \hbox{as}\quad \varepsilon\to 0.
\end{array}
\end{equation}
\par

Therefore, from (\ref{E-v}), (\ref{Eq-4.1-1})-(\ref{Eq-4.1-3}),
(\ref{Eq-4.1-6})-(\ref{Eq-4.1-12}), as $\varepsilon\to 0$, we have
$$
\begin{array}{l}
\int_0^T\int_{D}\nu u_{tt}(t,x) \phi(t,x)dxdt+ \int_0^{T}\int_{D}\nu
u_t(t,x) \phi(t,x)dxdt\\
+ \frac{1}{|Y|} \int_0^{T}\int_{D}\int_{Y^*} [\bigtriangledown_x
u(t,x)+\bigtriangledown_y u_1(t,x,y)][\bigtriangledown_x\phi(t,x)+
\bigtriangledown_y\Phi(t,x,y)] dydxdt \\
+ \int_0^{T}\int_{D}\nu
u(t,x)\phi(t,x)dxdt-\int_0^{T}\int_{D}\sin(\nu
u(t,x))\phi(t,x)dxdt\\
=\int_0^{T}\int_{D}\nu \varphi(t,x)dxdW_1(t),
\end{array}
$$
which implies that
\begin{equation}\label{HE}
\left\{
\begin{array}{l}
\nu u_{tt}(t,x)+ \nu u_t(t,x)-div_x\mathcal{A}(\bigtriangledown_x
u(t,x)) + \nu u(t,x)-\sin(\nu u(t,x)) =\nu \dot{W}_1,\\
\frac{\partial( \bigtriangledown_x u(t,x)+\bigtriangledown_y
u_1(t,x,y))}{\partial {\bf m}}=0,\quad \hbox{on}\quad \partial
Y^*-\partial Y.
\end{array}
\right.
\end{equation}
where ${\bf m}$ is the unit exterior norm vector on $\partial
Y^*-\partial Y$ and
$$
\mathcal{A}(\bigtriangledown_x u(t,x))= \frac{1}{|Y|} \int_{Y^*}
[\bigtriangledown_x u(t,x)+\bigtriangledown_y u_1(t,x,y)]dy,
$$
with $u_1$ satisfying for any $\Psi\in H_0^1(D_T; H_{per}^1(Y))$,
\begin{equation}\label{E-u1}
\left\{
\begin{array}{l}
\int_{Y^*} [\bigtriangledown_x u(t,x)+\bigtriangledown_y u_1(t,x,y)]
\bigtriangledown_y \Psi dy=0,\\
u_1\;\hbox{is}\;Y-\hbox{periodic}.
\end{array}
\right.
\end{equation}
Especially notice that Equation (\ref{E-u1}) has a unique solution
for any given $u(t,x)$, which implies that
$\mathcal{A}(\bigtriangledown_x u(t,x))$ is well-defined. Please
refer to \cite{FM} about the further properties of
$\mathcal{A}(\bigtriangledown_x u(t,x))$. Furthermore, from the
classic theory of stochastic partial differential equation, the
problem (\ref{HE}) is well-posed.
\par
In addition, from the classical homogenization theory (see \cite{CD,
FM}), we have
$$
u_1(t,x, y)=\sum\limits_{i=1}^3\frac{\partial u(t,x)}{\partial
x_i}(w_i(y)-e_i(y))
$$
where $\{e_i\}_{i=1}^3$ is the canonical basis of $\mathbb{R}^3$ and
$w_i$ is the solution of the following elementary cell problem
$$
\left\{
\begin{array}{l}
\bigtriangleup_y w_i(y)=0,\quad \hbox{in}\quad Y^*,\\
w_i-e_iy\quad \hbox{is}\quad Y-\hbox{periodic},\\
\frac{\partial w_i}{{\partial\bf n}}=0\quad \hbox{on}\quad \partial
S.
\end{array}
\right.
$$
Then $\mathcal{A}\bigtriangledown u=\mathcal{A}^*\bigtriangledown u$,
with $\mathcal{A}^*=(A_{ij}^*)$ being the classical homogenized
matrix defined as
\begin{equation}\label{Eq-A}
A_{ij}^*=\frac{1}{|Y|}\int_{Y^*}w_i(y)w_j(y)dy.
\end{equation}
\par
Then define $V(t,x)=\nu u(t,x)$.  Combining ({\ref{uv}) and
({\ref{HE}), we know that Theorem 4.1 holds. The proof is thus complete.
\hfill$\blacksquare$
\par
{\bf Remark 4.1}\quad {\it By the classic stochastic partial
differential equation theory, Equation (\ref{Eq-E}) is well-posed.
Here, we omit its proof. }

\end{document}